\documentclass[MSNbibl,number,citesort,dvips]{arxbj}
\usepackage{upgreek}
\usepackage{dcolumn}
\usepackage[left]{vmkcol}
\usepackage{graphicx}

\aid{0}
\volume{20}
\issue{2}
\pubyear{2014}
\firstpage{919}
\lastpage{957}
\doi{10.3150/13-BEJ510} 

\makeatletter
\newcolumntype{d}[1]{D{.}{.}{#1}}

\newcommand{\RMo}{\mathrm{o}}
\newcommand{\RMO}{\mathrm{O}}

\newcommand{\RMe}{\mathrm{e}}

\newcommand{\rrvert}{\vert}
\newcommand{\llvert}{\vert}

\newcommand{\dd}{\,\mathrm{d}}
\newcommand{\ddd}{\mathrm{d}}

\newcommand{\R}{\mathbb{R}}
\newcommand{\N}{\mathbb{N}}
\newcommand{\F}{\mathcal{F}}
\newcommand{\1}{\mathbf{1}}
\newcommand{\di}{\Delta_i}
\newcommand{\dn}{\Delta_n}
\newcommand{\ft}{\mathcal{F}}

\newtheorem{theorem}{Theorem}[section]
\newtheorem{proposition}[theorem]{Proposition}
\newtheorem{lemma}[theorem]{Lemma}

\newproclaim{assumption}[theorem]{Assumption}

\newremark{remark}[theorem]{Remark}
\newremark{example}[theorem]{Example}

\makeatother

\begin{document}
\begin{frontmatter}

\title{Efficient maximum likelihood estimation for L\'evy-driven
Ornstein--Uhlenbeck processes}
\runtitle{Efficient MLE for Ornstein--Uhlenbeck type processes}

\begin{aug}
\author{\fnms{Hilmar} \snm{Mai}\corref{}\ead[label=e1]{mai@wias-berlin.de}}
\runauthor{H. Mai} 
\address{Weierstrass Institute, Mohrenstr. 39, 10117 Berlin,
Germany. \printead{e1}}
\end{aug}

\received{\smonth{7} \syear{2012}}
\revised{\smonth{12} \syear{2012}}

%
\begin{abstract}
We consider the problem of efficient estimation of the drift parameter
of an Ornstein--Uhlenbeck type process driven by a L\'evy process when
high-frequency observations are given. The estimator is constructed
from the time-continuous likelihood function that leads to an explicit
maximum likelihood estimator and requires knowledge of the continuous
martingale part. We use a thresholding technique to approximate the
continuous part of the process. Under suitable conditions, we prove
asymptotic normality and efficiency in the H\'ajek--Le Cam sense for
the resulting drift estimator. Finally, we investigate the finite
sample behavior of the method and compare our approach to least squares
estimation.
\end{abstract}

%

\end{frontmatter}

\section{Introduction}

Let $(L_t,t \geq0)$ be a L\'{e}vy process on a filtered probability space
$(\Omega, \F, (\F_t), P)$ adapted to the filtration $(\F_t)_{t\geq
0}$. Denote by $(b,\sigma^2,\mu)$ the L\'evy--Khintchine triplet of
$L$. We call for every $a \in\R$ a strong solution $X$ to the
stochastic differential equation
%
\begin{equation}
\label{OUeq} \ddd X_t = -a X_t \dd t + \ddd L_t,\qquad t
\in\R_+,\qquad X_0 =\tilde{X},
\end{equation}
a L\'evy-driven Ornstein--Uhlenbeck (OU) process or Ornstein--Uhlenbeck
type process. The initial condition $\tilde X$ is assumed to be
independent of $L$. We consider the problem of estimating the mean
reversion parameter $a$ when observations $X_{t_1},\ldots,X_{t_n}$ on
an interval $[0,T_n]$ are given. It is well known that the drift of $X$
is identifiable only in the limit $T_n \to\infty$, even when
time-continuous observations are given. Therefore, we work under the
asymptotic scheme $T_n \to\infty$ and $\Delta_n = \max_{1 \leq i
\leq n-1} \{ |t_{i+1}- t_i|\} \downarrow0$ as $n \to\infty$.

The OU process serves us here as a toy model to understand the
interplay of jumps and continuous component of $X$ in this estimation
problem. This interplay is fundamental also for drift estimation in
more general models (cf. \cite{Mai2012}).

Ornstein--Uhlenbeck type processes have important applications in
various fields. In mathematical finance they are well know as a main
building block of the Barndorff--Nielsen--Shephard stochastic volatility
model (cf. \cite{BNS01}). But also in neuroscience they are popular
for the description of the membrane potential of a neuron (cf.
\cite{Jahn2011} and~\cite{Lansky1987}).

Estimation of L\'evy-driven Ornstein--Uhlenbeck processes has been
considered by several authors (see \cite{Masuda2010} and the
references therein) mostly when the driving L\'evy process is a
subordinator. Some examples are \cite{Jongbloed2005} on nonparametric
estimation of the L\'evy density of $L$, in \cite{Brockwell07} the
Davis--McCormick estimator was applied in the OU context and parametric
estimation based on a cumulant $M$-estimator was studied in
\cite{Jongbloed2006}. In \cite{Hu2009}, least squares estimation of the
drift parameter for an $\alpha$-stable driver is discussed, when no
Gaussian component is present. Reference \cite{Masuda2010} found that the rate
of convergence of the least absolute deviation estimator is either the
standard parametric rate, when $L$ has a Gaussian component, or is
faster than the standard rate, when $L$ is a pure jump process and
depends on the activity of the jumps. In \cite{Shimizu2006}, joint
parametric estimation of the drift and the L\'{e}vy measure was treated
via estimating
functions. Unfortunately, none of these methods lead to an efficient
estimator of the drift when $L$ is a general L\'evy
process.\looseness=1

To construct an efficient estimator, our starting point will be the
continuous time likelihood function. From this likelihood function, an
explicit maximum likelihood estimator can be derived, which is
efficient in the sense of the H\'ajek--Le Cam convolution theorem. In
the likelihood function the continuous martingale part of $X$ under the
dominating measure appears, which is not directly observed in our
setting. For discrete observations, we approximate the continuous part
of $X$ by neglecting increments that are larger than a certain
threshold that has to be chosen appropriately. We will call this
thresholding technique a jump filter. For this discretized likelihood
estimator with jump filtering, we prove asymptotic normality and
efficiency by showing that it attains the same asymptotic distribution
as the benchmark estimator based on time-continuous
observations.\looseness=1

This leads to the main mathematical question underlying this estimation
problem. Can we recover the continuous part of $X$ in the
high-frequency limit via jump filtering? If $L$ has only compound
Poisson jumps it is intuitively clear that the answer is yes. But when
$L$ has infinitely many small jumps in every finite interval this is a
much more challenging question. It turns out that even in this
situation jump filtering works under mild assumptions on the behavior
of the L\'evy measure around zero. The main condition here is that the
Blumenthal--Getoor index of the jump part is strictly less than two,
which is apparently a necessary condition in this context. In this
setting, jump filtering becomes possible since on a small time scale
increments of the continuous part and of the jump part exhibit a
different order of magnitude such that they can be distinguished via
thresholding. To control the small jumps of $L$ under thresholding, we
derive an estimate for the Markov generator of a thresholded pure jump
L\'evy process. Estimates of this type without thresholding were given,
among others, in \cite{Figueroa-Lopez2008} and \cite{Jacod2007}.

The problem of separation between continuous and jump part of a process
appears naturally in many situations. For example, in estimation of the
integrated volatility of a jump diffusion process via realized
volatility the quadratic variation of the jump component has to be
removed. This problem has been solved by thresholding in
\cite{Mancini2009} for Poisson jumps and in \cite{Mancini2011} for more
general jump behavior. Efficiency questions in this context in a simple
parametric model have been addressed in \cite{ASJ2007}. When
properties of the jump component are of interest thresholding
techniques are equally useful as demonstrated, among others, in \cite
{Ait-Sahalia2012}. In contrast to our discussion all these references
consider the separation problem for a finite and fixed observation
horizon $T_n = T < \infty$.

We also demonstrate in a simulation example that jump filtering leads
to a major improvement of the drift estimate for finite sample size.
Let us also mention that implementation of the drift estimator is
straightforward and computation time is not an issue even for large
data sets.

The paper is organized as follows: in Section \ref{sec2} we derive the
maximum likelihood estimator based on time-continuous observations,
give its asymptotic properties and obtain the efficient limiting
distribution for this estimation problem. Section \ref{chp5} deals
with estimating the drift parameter from discrete observations when $L$
has finite jump activity. In Section \ref{chp6} we build on the
results from Sections \ref{sec2} and \ref{chp5} to prove efficiency
also for possibly infinite jump activity. The finite sample behavior of
the estimator is investigated in Section \ref{chp7} based on simulated
data together with an analysis of the impact of the jump filter on the
performance of the estimator.

\section{Maximum likelihood estimation}\label{sec2}

Let us summarize some important facts on Ornstein--Uhlenbeck type
processes. It follows from It\^o's formula that an explicit solution of
(\ref{OUeq}) is given by
%
\begin{equation}
\label{OUsol} X_t = \RMe^{-a t} X_0 + \int
_0^t \RMe^{-a (t-s)} \dd L_s,\qquad t \in
\R_+.
\end{equation}
The integral in (\ref{OUsol}) can by partial integration be defined
path-wise as a Riemann--Stieltjes integral, since the integrand is of
finite variation (see \cite{Dudley}, e.g.). This solution to
equation (\ref{OUeq}) is unique up to indistinguishability.

When $L$ is non-deterministic it was shown in \cite{Jacod1985} that
equation (\ref{OUeq}) admits a causal stationary solution (cf. also
\cite{Masuda2007} and \cite{SY}) if and only if
%
\begin{equation}
\label{OUstatcond} \int_{|x| > 1} {\log}|x| \mu(\ddd x) < \infty
\quad\mbox{and}\quad a > 0.
\end{equation}
Under these conditions $X$ has a unique invariant distribution $G$ and
$X_t \stackrel{\mathcal{D}}{\longrightarrow}X_\infty\sim G$ as $t
\to\infty$.

The Ornstein--Uhlenbeck process $X$ exhibits a modification with c\`adl\`ag
paths and hence it induces a measure $P^a$ on the space $D[0,\infty)$
of c\`adl\`ag functions on the interval $[0,\infty)$. Denote by $P_t^a$
the restriction of $P^a$ to the $\sigma$-field $\ft_t$. If $\sigma^2
>0$ then these induced measure are locally equivalent (cf.
\cite{Sorensen1991}) and the corresponding Radon--Nikodym derivative or
likelihood function is given by
\[
\frac{\ddd P_t^a}{\ddd P_t^0} = \exp\biggl( -\frac{ a}{ \sigma^2} \int_0^T
X_{s} \dd X_s^c - \frac{a^2}{2 \sigma^2} \int
_0^T X_{s}^2 \dd s
\biggr),
\]
where $X^c$ denotes the continuous $P^0$-martingale part of $X$. This
leads to the explicit maximum likelihood estimator
%
\begin{equation}
\label{MLE}
\hat{a}_T = -\frac{\int_0^T X_{s} \dd X_s^c}{\int_0^T X_{s}^2 \dd s}
\end{equation}
for $a$ when the process is fully observed on $[0,T]$. The estimator
$\hat a_T$ cannot be applied in this form, since time-continuous
observations are usually not available in most applications. Therefore,
we will develop in the next section a discrete version of $\hat a_T$
and prove its efficiency. The main challenge there will be that the
continuous part $X^c$ is not directly observed and hence has to be
approximated from discrete observations of $X$ via jump filtering.

However, for us $\hat a_T$ will serve as a benchmark for the estimation
problem with discrete observations. Asymptotic normality and efficiency
in the H\'ajek--Le Cam sense of $\hat a_T$ follow easily from general
results for exponential families of stochastic processes (cf.
\cite{KS97} and \cite{Mai2012}). These results provide an efficiency bound
for the case of discrete observations in the next section. Let us
summarize them in the following theorem.
%
\begin{theorem} \label{consistencyThm}
\textup{(i)} Under the condition $\sigma^2 > 0$ the estimator $\hat
{a}_T$ exists uniquely and is strongly consistent under $P^a$, that is,
\[
\hat{a}_T \stackrel{a.s.} {\longrightarrow}a
\]
under $P^a$ as $T \to\infty$.

\mbox{}\hphantom{\textup{i}}\textup{(ii)} Suppose that additionally (\ref{OUstatcond}) holds and
that $X$ has bounded second moments such that the invariant
distribution satisfies $E_a [X_\infty^2] < \infty$. Then under $P^a$
\[
\sqrt{T} (\hat{a}_T -a) \stackrel{\mathcal{D}} {\longrightarrow}N
\biggl(0, \frac{\sigma
^2}{E_a[X_\infty^2]} \biggr)
\]
and
%
\begin{equation}
\label{oustdnormal} \sigma^{-1} S_T^{1/2} (\hat
a_T - a) \stackrel{\mathcal{D}} {\longrightarrow}N(0,1)
\end{equation}
as $T \to\infty$, where $S_T = \int_0^T X_s^2 \dd s $.

\textup{(iii)} The statistical experiment $\{P^a, a \in\R_+ \}$ is
locally asymptotically normal.

\mbox{}\hspace*{0.7pt}\textup{(iv)} The estimator $\hat a_T$ is asymptotically efficient in
the sense of H\'ajek--Le Cam.
\end{theorem}
For a proof we refer to \cite{Mai2012}, Section 4.2.
%
\begin{remark}
When $\sigma^ 2$ is known or a consistent estimator is at hand, we can
use (\ref{oustdnormal}) to construct confidence intervals for $a$.
\end{remark}

\section{Discrete observations: Finite activity} \label{chp5}

In this section, we consider the estimation of $a$ for discrete
observations. The maximum likelihood estimator for the drift given in
(\ref{MLE}) involves the continuous martingale part that is unknown
when only discrete observations are given. Hence, we will approximate
the continuous part of the process by removing observations that most
likely contain jumps. We restrict our attention in this section to the
case that the driving L\'{e}vy process has jumps of finite activity. The
jump filtering technique provides us in the high-frequency limit an
asymptotically normal and efficient estimator. Based on these results,
we will treat the general case of an infinitely active jump component
in Section \ref{chp6}.

\subsection{Estimator and observation scheme} \label{oudsec}

Let $X$ be an Ornstein--Uhlenbeck process defined by (\ref{OUeq}) and
suppose we observe $X$ at discrete time points $0 = t_1 < t_2 < \cdots
< t_{n} =T_n$ such that $T_n \to\infty$ as well as $\Delta_n = \max_{1
\leq i \leq n-1} \{ |t_{i+1}- t_i|\} \downarrow0$ and $n \dn
T_n^{-1} = \RMO(1)$ as $n \to\infty$. The last condition assures that
the number of observations $n$ does not grow faster than $T_n \dn
^{-1}$. It can always be fulfilled by neglecting observations and will
simplify the formulation of the proof considerably. Denote by
$(b,\sigma^2, \mu)$ the L\'evy--Khintchine triplet of $L$. Assume
throughout this section that $\lambda= \mu(\R) < \infty$ for the L\'evy measure $\mu$.

By deleting increments that are larger than a threshold $v_n > 0$ we
filter increments that most likely contain jumps and thus approximate
the continuous part with the remaining increments. Applied to the
time-continuous likelihood estimator (\ref{MLE}) this method leads to
the following estimator:
%
\begin{equation}
\label{baranchp5} \bar{a}_n:= -\frac{\sum_{i=0}^{n-1}X_{t_i} \Delta_i X
\1_{\{
|\Delta_i X|
\leq v_n\}}}{\sum_{i=0}^{n-1}X_{t_i}^2 (t_{i+1} - t_i)}.
\end{equation}
Here $v_n >0, n\in\N$, is a cut-off sequence that will be chosen as a
function of the maximal distance between observations $\dn$ and $\di X
= X_{t_{i+1}} - X_{t_i}$.

In the finite activity case, the jump part $J$ of $L$ can be written as
a compound Poisson process
\[
J_t = \sum_{i= 1}^{N_t}
Z_i,
\]
where $N$ is a Poisson process with intensity $\lambda$ and the jump
heights $Z_1, Z_2, \ldots$ are i.i.d. with distribution $F$.

\subsection{Asymptotic normality and efficiency}

The indicator function that appears in $\bar a_n$ deletes increments
that are larger than $v_n$. In \cite{Mancini2009}, it was shown that
increments of the continuous part of $X$ over an interval of length
$\dn$ are with high probability smaller than $\dn^{1/2}$. Hence, we
set $v_n = \dn^\beta$ for $\beta\in(0,1/2)$ to keep the continuous
part\vspace*{2pt} in the limit unaffected by the threshold. In order to be able to
choose $v_n$ such that $X^c_n = \sum_{i=0}^{n-1}\di X \1_{\{|\Delta_i
X| \leq
v_n\}}$ approximates the continuous martingale part in the limit, we
make the following assumptions on the jumps of $L$ and the observation scheme.
%
\begin{assumption} \label{ajumps}
\textup{(i)} Suppose that $\mu$ and $a$ satisfy (\ref{OUstatcond}),
the drift $b=0$, the process $X$ has bounded second moments,

\mbox{}\hphantom{\textup{i}}\textup{(ii)} the distribution $F$ of the jump heights is such that
\[
F\bigl(\bigl(-2 \Delta_n^\beta, 2 \Delta_n^\beta
\bigr)\bigr) = \mathrm{o}\bigl(T_n^{-1}\bigr),
\]

\textup{(iii)} and there exists $\beta\in(0,1/2)$ such that the
maximal distance between observations satisfies $T_n \dn^{(1-2 \beta)
\wedge({1/2}) } = \RMo(1)$.
\end{assumption}
%
\begin{remark} \label{jrem1}
Suppose that $F$ has a bounded Lebesgue density $f$. Then
\[
F\bigl(\bigl( -\Delta_n^\beta, \Delta_n^\beta
\bigr)\bigr) = \mathrm{O}\bigl(\dn^\beta\bigr)
\]
and Assumption \ref{ajumps}(ii) becomes $\dn^\beta T_n =\RMo (1)$.
Together with Assumption \ref{ajumps}(i), we obtain that $\beta= 1/3$ leads to
an optimal compromise between Assumption \ref{ajumps}(i) and (ii).
\end{remark}
%
\begin{remark}
Assumption \ref{ajumps}(iii) means here that for given $T_n \to\infty
$ we require $\dn\downarrow0$ fast enough such that there exists
$\beta\in(0,1/2)$:
\[
T_n \dn^{1- 2 \beta} = \RMo (1) \quad\mbox{and}\quad T_n
\dn^{1/2} = \RMo(1).
\]
Of course one of these two conditions will be dominating and determine
the order of $\dn$.
\end{remark}
%
\begin{remark}
Assumption \ref{ajumps}(ii) gives a lower bound for the choice of the
threshold $\beta$. At the same time Assumption \ref{ajumps}(iii) limits the
range of possible $\beta$'s from above, since the available frequency
of observations, that is, the order of $\dn$, may be limited in specific
applications. Hence, the distribution $F$, the observation length $T_n$
and frequency $\dn$ fix a range for the choice of $\beta$. At this
point the question of a data driven method to choose $\beta$ arises,
but this will not be considered in this work. The condition $b=0$ is
necessary in this context, since otherwise there is no hope to recover
the continuous martingale part via jump filtering.
\end{remark}

The following theorem gives as the main result of this section a
central limit theorem for the discretized MLE with jump filter.
%
\begin{theorem}\label{discTCLT}
Suppose that Assumption \ref{ajumps} holds and $\sigma^2 > 0$. Set
$v_n = \Delta_n^\beta$ for $\beta\in(0,1/2)$, then
\[
T_n^{1/2} (\bar{a}_n - a) \stackrel{
\mathcal{D}} {\longrightarrow} N \biggl(0, \frac{\sigma^2}{E_a[X_\infty
^2]} \biggr) \qquad\mbox{as $n
\to\infty$}.
\]
The estimator $\bar a_n$ is asymptotically efficient.
\end{theorem}
%
\begin{remark}
Asymptotic efficiency follows immediately from Theorem \ref
{consistencyThm} and the first statement of Theorem \ref{discTCLT}.
\end{remark}

\subsection{Proofs}
We divide the proof of the theorem into several lemmas. First of all,
we need a probability bound for the event that the continuous component
of $X$ exceeds a certain threshold.

By the L\'evy--It\^o decomposition (cf. \cite{Sato99}) and since $b=0$
in our setting the driving L\'{e}vy process can be decomposed as $L = W +
J$, where $W$ is a standard Wiener process and $J$ is a pure jump L\'{e}vy
process independent of $W$. Denote by $D$ the drift component of $X$,
that is,
\[
D_t = -a\int_0^t X_s
\dd s.
\]

\begin{lemma} \label{pbound}
Let $\sup_{s\geq0}\{ E[|X_s|^l]\}< \infty$ for some $l \geq1$. For
any $\delta\in(0,1/2)$ and $i \in\{ 1,\ldots, n-1\}$, we have
\[
P\bigl(|\di W + \di D| > \Delta_n^{1/2- \delta}\bigr) =\mathrm{O} \bigl(
\Delta_n^{l(1/2 + \delta)} \bigr) \qquad\mbox{as $n \to\infty$}.
\]
\end{lemma}
\begin{pf}
In the first step, we separate $\di W$ and $\di D$:
\[
P \bigl(|\di W + \di D| > \Delta_n^{1/2- \delta} \bigr) \leq P
\biggl(|\di W| >\frac{ \Delta_n^{1/2 - \delta}}{2} \biggr) + P \biggl
(|\di D| >\frac{ k \Delta_n^{1/2}}{2}
\biggr).
\]
By Lemma 22.2 in \cite{Klenke},
\[
P \biggl(|\di W| >\frac{ \Delta_n^{1/2 - \delta}}{2} \biggr) \leq2 \dn
^\delta
\RMe^{-{1}/({8 \dn^\delta})}.
\]
It follows from Jensen's inequality that
\[
\biggl\llvert\int_{t_i}^{t_{i+1}} X_s \dd
s \biggr\rrvert^l \leq\dn^{l-1} \int_{t_i}^{t_{i+1}}
|X_s|^l \dd s.
\]
This leads to
\[
E\bigl[|\di D|^l\bigr] \leq\dn^{l-1} a^l \int
_{t_i}^{t_{i+1}} E\bigl[|X_s|^l
\bigr] \dd s \leq\dn^l a^l \sup_{s \in[t_i, t_{i+1}]}
\bigl\{ E\bigl[|X_s|^l\bigr]\bigr\}.
\]
Finally, Markov's inequality yields
\[
P \biggl(|\di D| >\frac{ \Delta_n^{1/2- \delta}}{2} \biggr) \leq a^l
\sup
_{s \in[t_i, t_{i+1}]}\bigl\{ E\bigl[|X_s|^l\bigr]
\bigr\} \frac{2^l \Delta
_n^{l}}{ \dn^{l (1/2 - \delta)}} = \mathrm{O} \bigl( \Delta_n^{l(1/2 +
\delta)} \bigr).
\]
\upqed\end{pf}

\subsubsection{Jump filtering}
First, we will investigate how to choose the cut-off sequence $v_n$ in
order to filter the jumps. Define for $n \in\N$ and $i \in\{1,\ldots, n\}$ the following sequence of events
\[
A_n^i = \bigl\{ \omega\in\Omega\dvt  \1_{\{ |\Delta_i X| \leq v_n \}} (
\omega) = \1_{\{\Delta_i N =0 \}} (\omega) \bigr\}.
\]
Here $N$ denotes the counting measure that counts the jumps of $L$.

\begin{lemma}\label{Ain}
Suppose that Assumption \ref{ajumps} holds and set $v_n = \Delta
_n^\beta, \beta\in(0,1/2)$, then it follows that for $A_n = \bigcap
_{i=1}^{n} A_n^i$ we have
\[
P(A_n) \to1 \qquad\mbox{as $n \to\infty$.}
\]
\end{lemma}

\begin{pf}
Observe that
\[
P\bigl(A_n^c\bigr) = P \Biggl( \bigcup
_ {i = 1}^n\bigl(A_n^i
\bigr)^c \Biggr) \leq\sum_{i
= 1}^n
P\bigl(\bigl(A_n^i\bigr)^c\bigr).
\]
By setting
\begin{eqnarray*}
K_n^i &=& {\bigl\{ |\Delta_i X| \leq
v_n \bigr\}},
\\
M_n^i &=& {\{\Delta_i N =0 \}},
\end{eqnarray*}
we can rewrite $(A_n^i)^c$ as
\[
\bigl(A_n^i\bigr)^c = \{
\1_{ K_n^i } \neq\1_{ M_n^i } \} = \bigl(K_n^i
\setminus M_n^i\bigr) \cup\bigl(M_n^i
\setminus K_n^i\bigr).
\]
Here the events $K_n^i \setminus M_n^i$ and $M_n^i \setminus K_n^i$
correspond to the two types of errors that can occur when we search for
jumps. In the first case, we miss a jump and in the second case we
neglect an increment although it does not contain any jumps. Next, we
are going to bound the probability of both errors:
%
\begin{equation}
\label{Aidecomp} P\bigl(\bigl(A_n^i
\bigr)^c\bigr) = P\bigl(K_n^i \setminus
M_n^i\bigr) + P\bigl(M_n^i
\setminus K_n^i\bigr).
\end{equation}
Set $\di= t_{i+1} - t_i$. For the first type of error, we obtain
%
\begin{eqnarray}
\label{jdecomplem} P\bigl(K_n^i \setminus
M_n^i\bigr) &=& P\bigl(|\Delta_i X| \leq
v_n, \Delta_i N >0\bigr)
\nonumber
\\
&=& \sum_{j=1}^\infty \RMe^{-\lambda\di}
\frac{(\lambda
\di)^j}{j!} P\bigl(|\Delta_i X| \leq v_n |
\Delta_i N =j\bigr)
\\
&\leq& P(\di N =1) P \bigl(|\Delta_i X| \leq v_n |
\Delta_i N =1\bigr) + \mathrm{O}\bigl(\Delta_n^2\bigr)\nonumber
\end{eqnarray}
and
%
\begin{eqnarray}\label{in1}
P\bigl(|\di X| \leq v_n| \di N =1\bigr) &\leq& P\bigl(|\di X| \leq
v_n, |\di J| > 2v_n | \di N =1\bigr)
\nonumber\\[-8pt]\\[-8pt]
&&{}+ P\bigl(|\di X| \leq v_n, |\di J| \leq2v_n | \di N =1\bigr).
\nonumber
\end{eqnarray}
The first term on the right-hand side is bounded by
%
\begin{eqnarray}
\label{onejump1} &&
P\bigl(|\di X| \leq v_n, |\di J| > 2v_n |
\di N =1 \bigr)
\nonumber
\\
&&\quad= P \bigl(|\Delta_i W + \Delta_i J + \Delta_i
D| \leq v_n, |\Delta_i J|> 2v_n | \di N =
1\bigr)
\nonumber\\[-8pt]\\[-8pt]
&&\quad\leq P \bigl(|\Delta_i W + \Delta_i D| > v_n,
\di N =1\bigr) P(\di N =1)^{-1}
\nonumber
\\
&&\quad\leq P \bigl(|\Delta_i W + \Delta_i D| > v_n\bigr)
P(\di N =1)^{-1} = P(\di N =1)^{-1} \mathrm{O}\bigl(
\Delta_n^{2-2 \beta}\bigr),\nonumber
\end{eqnarray}
where we used Lemma \ref{pbound} with $l=2$.
Denote by $F$ the distribution of the jump heights of $J$. Then we
obtain for the second term on the right-hand side of (\ref{in1})
\[
P\bigl(|\di X| \leq v_n, |\di J| \leq2v_n | \di N =1\bigr) \leq
P\bigl(|\di J| \leq2v_n | \di N =1\bigr) = F\bigl((-2v_n,2v_n)
\bigr).
\]
For the second addend in (\ref{Aidecomp}) it follows by independence
of $W$ and $J$ that
\begin{eqnarray*}
P\bigl(M_n^i \setminus K_n^i
\bigr) &=& P\bigl(|\Delta_i X| > v_n, \Delta_i N
=0\bigr)
\\
&\leq& P\bigl(|\Delta_i W + \Delta_i D| > v_n\bigr).
\end{eqnarray*}
Lemma \ref{pbound} yields
%
\begin{equation}
\label{onejump2} P\bigl(|\Delta_i W + \Delta_i D| >
v_n\bigr) = \mathrm{O}\bigl(\Delta_n^{2-2 \beta}\bigr).
\end{equation}
Finally, (\ref{jdecomplem}), (\ref{onejump1}) and (\ref{onejump2})
lead to
\[
P\bigl(\bigl(A_n^i\bigr)^c\bigr) \leq F
\bigl(\bigl(-2 \Delta_n^{\beta}, 2 \Delta_n^{\beta}
\bigr)\bigr) \Delta_n + \mathrm{O}\bigl(\Delta_n^{2-2 \beta}
\bigr)
\]
such that the statement follows, since we have shown that
\[
P\bigl(A_n^c\bigr) \leq\sum
_{i = 1}^n P\bigl(\bigl(A_n^i
\bigr)^c\bigr) \leq \RMO(T_n) F\bigl(\bigl(-2 \Delta
_n^\beta, 2\Delta_n^\beta\bigr)\bigr)
+ \mathrm{O}\bigl(T_n \dn^{1-2 \beta}\bigr).
\]
\upqed\end{pf}

\subsubsection{Approximation of the continuous martingale part}
%
\begin{lemma} \label{JApp}
Under Assumption \ref{ajumps}, we obtain
\[
\Biggl\llvert\sum_{i=0}^{n-1}X_{t_i}
\bigl(\Delta_i X \1_{\{|\Delta_i X| \leq
v_n\}} - \Delta_i
X^c\bigr) \Biggr\rrvert= \mathrm{O}_p \bigl(T_n
\dn^{1/2}\bigr)
\]
as $n \to\infty$.
\end{lemma}

\begin{pf}
On $A_n$ from Lemma \ref{Ain}, we have
%
\begin{equation}
\label{incdiff} \sum_{i=0}^{n-1}X_{t_i}
\bigl( \Delta_i X \1_{\{|\Delta_i X| \leq
v_n\}} - \Delta_i
X^c \bigr) = \sum_{i=0}^{n-1}X_{t_i}
\bigl( \Delta_i X \1 _{\{
\Delta_i N = 0\}} - \Delta_i
X^c \bigr).
\end{equation}
By Lemma \ref{Ain}, we have $P(A_n) \to1$ as $n \to\infty$. Observe
now that the difference of the increments on the right-hand side of
(\ref{incdiff}) is unequal to zero only if a jump occurred in that
interval, that is,
\[
\Delta_i X \1_{\{\Delta_i N = 0\}} - \Delta_i
X^c = \cases{ - \Delta_i X^c; &\quad $
\Delta_i N > 0$,
\cr
0; &\quad $\Delta_i N =0$.}
\]
Define $C_i^n = \{\Delta_i N > 0\}$ and observe that
\[
E \Biggl\llvert\1_{A_n} \sum_{i=0}^{n-1}X_{t_i}
\bigl(\Delta_i X \1_{\{\Delta
_i N = 0\}} - \Delta_i
X^c\bigr) \Biggr\rrvert= E \Biggl\llvert\sum
_{i=0}^{n-1}X_{t_i} \Delta_i
X^c \1_{A_n \cap
C_i^n} \Biggr\rrvert.
\]
The $i$th increment of $X^c$ can be written as $\Delta_i X^c = \Delta
_i W + \Delta_i D$. Therefore,
\begin{eqnarray*}
E \Biggl\llvert\sum_{i=0}^{n-1}X_{t_i}
\Delta_i X^c \1_{A_n \cap C_i^n} \Biggr\rrvert&\leq&\sum
_{i=0}^{n-1}E \bigl[ \bigl|X_{t_i} (
\Delta_i W + \Delta_i D)\bigr| \1_{A_n
\cap C_i^n} \bigr]
\\
&\leq&\sum_{i=0}^{n-1}E \bigl[
\bigl(|X_{t_i} \Delta_i W| + |X_{t_i}
\Delta_i D|\bigr) \1_{ C_i^n} \bigr].
\end{eqnarray*}
The number of jumps of $J$ follows a Poisson process with intensity
$\lambda$ such that $P(C_i^n) \leq\Delta_n \lambda$. The
independence of $N \perp W$ and $\di N \perp X_{t_i}$ yields
\[
\sum_{i=0}^{n-1}E \bigl[ |X_{t_i}
\Delta_i W| \1_{ C_i^n} \bigr] = \sum_{i=0}^{n-1}
E \bigl[ |X_{t_i}|\bigr] E\bigl[| \Delta_i W|\bigr] P
\bigl(C_i^n\bigr) \leq \mathrm{O}\bigl(T_n
\Delta_n^{1/2}\bigr).
\]
Finally, by H\"older's inequality
\[
\sum_{i=0}^{n-1}E \bigl[ |X_{t_i}
\Delta_i D| \1_{ C_i^n} \bigr] \leq\sum
_{i=0}^{n-1}E \bigl[ X_{t_i}^2(
\Delta_i D)^2 \bigr]^{1/2} P
\bigl(C_i^n\bigr)^{1/2} = \mathrm{O}
\bigl(T_n \Delta_n^{1/2}\bigr).
\]
\upqed\end{pf}

\subsubsection{Central limit theorem for the discretized estimator}
To prove Theorem \ref{discTCLT}, we show next that when we discretize
the time-continuous estimator $\hat a_T$ as
\[
\hat{a}_n = - \frac{\sum_{i=1}^n X_{t_i} \Delta_i X^c}{\sum_{i=1}^n X_{t_i}^2
(t_{i+1} - t_i)},
\]
then $\hat a_n$ attains the same asymptotic distribution as $\hat a_T$
itself. In the last step, we will then show that the discretized MLE
and the estimator with jump filter show the same limiting
behavior.\vadjust{\goodbreak}
%
\begin{lemma} \label{discCLT}
If Assumption \ref{ajumps} is fulfilled, then the convergence
\[
T_n^{1/2}( \hat a_n - a ) \stackrel{
\mathcal{D}} {\longrightarrow}N \bigl( a, \sigma^2 E_a
\bigl[X_\infty^2\bigr]^{-1} \bigr) \qquad\mbox{as $n \to
\infty$}
\]
holds under $P^a$.
\end{lemma}
\begin{pf}
Let $W$ denotes a $P^a$-Wiener process. The continuous $P^0$-martingale
part can be written as
\[
X_t^c = \sigma W_t - a \int
_0^t X_s \dd s.
\]
This leads to the decomposition
\[
T_n^{1/2} (\breve a_n - a)=
T_n^{1/2} a \biggl(\frac{\sum_{i=0}^{n-1}X_{t_i}
\int_{t_i}^{t_{i+1}}X_s \dd s}{\sum_{i=0}^{n-1}X_{t_i}^2 \di} -1 \biggr
) -
T_n^{1/2} \frac{\sigma\sum_{i=0}^{n-1}X_{t_i} \di W}{\sum
_{i=0}^{n-1}X_{t_i}^2 \di} = S_n^1
- S_n^2.
\]
We will show now that $S_n^1 \stackrel{p}{\longrightarrow}0$ and
$S_n^2 \stackrel{\mathcal{D}}{\longrightarrow}N(0, \sigma^2
E_a[X_\infty^2]^{-1})$ as $n \to\infty$ such that the statement of
the proposition follows. Define $\lfloor t\rfloor_n = \max_{i \leq n} \{
t_i | t_i \leq t \}$. Let us first consider convergence of $S_n^1$.
Observe that
%
\begin{equation}
\label{S1nadecomp} S_n^1/a = \frac{T_n^{-1/2}(\sum_{i=0}^{n-1}X_{t_i}
\int_{t_i}^{t_{i+1}}X_s\dd s - \sum_{i=0}^{n-1}X_{t_i}^2 \di)}{T_n^{-1} \sum_{i=0}^{n-1}X_{t_i}^2
\di}.
\end{equation}
For the numerator, we obtain
%
\begin{eqnarray}
\label{enumesti} T_n^{-1/2} E \Biggl[ \Biggl\llvert\sum
_{i=0}^{n-1}X_{t_i} \int
_{t_i}^{t_{i+1}}X_s \dd s - \sum
_{i=0}^{n-1}X_{t_i}^2 \di\Biggr
\rrvert\Biggr] &\leq& T_n^{-1/2} \int_0^{T_n}
E_a \bigl[\bigl|X_{\lfloor t \rfloor_n} X_t - X_{\lfloor t
\rfloor_n}^2\bigr|
\bigr] \dd t\qquad
\nonumber\\[-8pt]\\[-8pt]
&=& \mathrm{O}\bigl(T_n^{1/2} \dn^{1/2}\bigr)\nonumber
\end{eqnarray}
such that the numerator converges to zero in $L^1$. A similar estimate
for the denominator yields
\[
T_n^{-1} E_a \Biggl[ \Biggl\llvert\int
_0^{T_n} X_t^2 \dd t - \sum
_{i=0}^{n-1} X_{t_i}^2
\di\Biggr\rrvert\Biggr] = \mathrm{O}\bigl(\Delta_n^{1/2}\bigr),
\]
and since the ergodic theorem implies that $T_n^{-1} \int_0^{T_n}
X_t^2 \dd t \stackrel{p}{\longrightarrow}E_a [X_\infty^2]$ as $n \to
\infty$, we conclude
%
\begin{equation}
\label{denomconv} T_n^{-1} \sum
_{i=0}^{n-1}X_{t_i}^2 \di
\stackrel{p} {\longrightarrow}E_a \bigl[X_\infty^2
\bigr]
\end{equation}
as $n \to\infty$. This convergence together with (\ref{S1nadecomp})
and the estimate (\ref{enumesti}) imply that $S_n^1 \stackrel
{p}{\longrightarrow}0$ as $n \to
\infty$.

It remains to prove convergence of $S_n^2$. From It\^o's isometry and
stationarity of $X$, we obtain for the numerator of $S_n^1$ that
\begin{eqnarray*}
T_n^{-1}E_a \Biggl[ \Biggl( \int
_0^{T_n} X_t \dd W_t - \sum
_{i=0}^{n-1}X_{t_i} \di W
\Biggr)^2 \Biggr] &=&T_n^{-1} E_a
\biggl[ \biggl( \int_0^{T_n} (X_t -
X_{\lfloor t \rfloor_n}) \dd W_t \biggr)^2 \biggr]
\\
&=& T_n^{-1} E_a \biggl[ \int
_0^{T_n} (X_t - X_{\lfloor t \rfloor
_n})^2
\dd t \biggr] \\
&=& T_n^{-1} \int_0^{T_n}
E_a \bigl[(X_t - X_{\lfloor t \rfloor_n})^2 \bigr]
\dd t = \RMO(\Delta_n).
\end{eqnarray*}
The numerator of $S_n^2$ is a continuous martingale and its quadratic
variation converges due to the ergodic theorem to the second moment of
$X$. The martingale central limit theorem implies now
\[
T_n^{-1/2} \sigma\int_0^{T_n}
X_t \dd W_t \stackrel{\mathcal{D}} {\longrightarrow}N
\bigl(0,\sigma^2 E_a\bigl(X_\infty^2
\bigr)\bigr)
\]
such that also
\[
T_n^{-1/2} \sigma\sum_{i=0}^{n-1}X_{t_i}
\di W \stackrel{\mathcal{D}} {\longrightarrow}N\bigl(0,\sigma^2
E_a\bigl(X_\infty^2\bigr)\bigr)
\]
as $n \to\infty$. This convergence together with (\ref{denomconv})
and Slutsky's lemma lead to
\[
S_n^1 \stackrel{\mathcal{D}} {\longrightarrow}N
\bigl(0, \sigma^2 E_a\bigl[X_\infty^2
\bigr]^{-1}\bigr)
\]
as $n \to\infty$. This completes the proof.
\end{pf}

\begin{pf*}{Proof of Theorem \ref{discTCLT}}
By Lemma \ref{discCLT} $T_n^{1/2} (\hat{a}_n - a) \stackrel{\mathcal
{D}}{\longrightarrow}N(0, \frac
{\sigma^2}{E_a[X_\infty^2]})$ as $n \to\infty$. By Slutsky's lemma,
it remains to show
%
\begin{equation}
\label{econv} T_n^{1/2} (\bar{a}_n -
\hat{a}_n) \stackrel{p} {\longrightarrow}0 \qquad\mbox{as $n \to\infty$.}
\end{equation}
Observe that
\[
T_n^{1/2} (\bar{a}_n - \hat{a}_n)
= T_n^{1/2} \biggl( \frac{\sum_{i=1}^n X_{t_i} \Delta_i X \1_{\{|\Delta
_i X| \leq v_n\}} - \sum_{i=1}^n X_{t_i} \Delta_i X^c}{\sum_{i=1}^n
X_{t_i}^2 \di} \biggr).
\]
By Lemma \ref{JApp} we obtain under $P_a$ that
\[
T_n^{-1/2} \Biggl(\sum_{i=1}^n
X_{t_i} \Delta_i X \1_{\{|\Delta_i X|
\leq v_n\}} - \sum
_{i=1}^n X_{t_i} \Delta_i
X^c \Biggr) \stackrel{p} {\longrightarrow}0 \qquad\mbox{as $n\to\infty$}
\]
and
\[
T_n^{-1} \sum_{i=1}^n
X_{t_i}^2 \di\stackrel{p} {\longrightarrow}E_a
\bigl[X_\infty^2\bigr],
\]
such that (\ref{econv}) follows.
\end{pf*}

\section{Discrete observations: Infinite activity} \label{chp6}

In this section, we generalize the results from Section \ref{chp5} to
the case that the jump part of the driving L\'{e}vy process can be of
infinite activity. We give conditions on the L\'{e}vy measure and suitable
rates for the cut-off sequence that ensure separation in the
high-frequency limit between jump part and continuous part. Under these
conditions, we will then prove asymptotic normality and efficiency of
the drift estimator $\bar a_n$ given in (\ref{baranchp5}).

The observation scheme considered here will be like in Section \ref{oudsec},
that is, $0 = t_1 < t_2 < \cdots< t_{n} =T_n$ such that $T_n \to
\infty$ as well as $\Delta_n = \max_{1 \leq i \leq n-1} \{ |t_{i+1}-
t_i|\} \downarrow0$ and $n \dn T_n^{-1} = \RMO(1)$ as $n \to\infty$.

\subsection{Asymptotic normality and efficiency}

In this section, we state as the main result of this paper a CLT for
the estimation error of $\bar a_n$. The limiting distribution will
imply asymptotic efficiency of $\bar a_n$. But before we can formulate
the theorem, we introduce some notation and mild assumptions on the
jump part of $L$ that enable us to separate the jump part and
continuous part via jump filtering.

Let $N$ denote the Poisson random measure associated to the jump part
of $L$. The jump component $J$ of $X$, the components $M$ of jumps
smaller than one and $U$ of jumps larger than one and the drift $D$ are
given by
%
\begin{eqnarray}
\label{jnotation} J_t &=& \int_0^t
\int_{-\infty}^\infty x \bigl(N(\ddd x,\ddd s)- \mu(\ddd x) \lambda
(\ddd s)\bigr),
\nonumber
\\
M_t &=& \int_0^t \int
_{-1}^1 x \bigl(N(\ddd x,\ddd s)- \mu(\ddd x) \lambda(\ddd s)
\bigr),
\nonumber\\[-8pt]\\[-8pt]
U_t &=& J_t - M_t,
\nonumber
\\
D_t &=& -a \int_0^t
X_s \dd s,
\nonumber
\end{eqnarray}
respectively. Owing to this decomposition of $X$ we can apply the
results from Chapter \ref{chp5} to $D$, $W$ and $U$ and thus can focus
on $M$. To control the small jumps of $M$, we impose the following
assumption on the L\'{e}vy measure $\mu$.
%
\begin{assumption} \label{aj}
\textup{(i)} Suppose that (\ref{OUstatcond}) holds, $b=0$ and $X$ has
bounded second moments.

\mbox{}\hphantom{\textup{i}}\textup{(ii)} There exists an $\alpha\in(0,2)$ such that as
$v\downarrow0$
%
\begin{equation}
\label{smallj} \int_{-v}^v x^2
\mu(\ddd x) = \mathrm{O}\bigl(v^{2-\alpha}\bigr).
\end{equation}

\textup{(iii)} There exists $\eta> 0$ and $n_0 \in\N$ such that for
all $\varepsilon\leq\eta$ and $n \geq n_0$
\[
E[\di M \1_{\{|\di M | \leq\varepsilon\}}] =0 \qquad\mbox{for all $i \in\{
1,\ldots, n-1\}$}.
\]
\end{assumption}

\begin{remark}
Assumption \ref{aj}(ii) controls the intensity of small jumps, which
is determined by the mass of $\mu$ around the origin. When $\gamma$
denotes the Blumenthal--Getoor index of $L$ defined by
\[
\gamma= \inf_{c \geq0} \biggl\{ \int_{|x| \leq1}
|x|^c \mu(\ddd x) < \infty\biggr\} \leq2,
\]
then $\alpha= \gamma$ satisfies (\ref{smallj}), that is, Assumption \ref
{aj}(i) states that the Blumenthal--Getoor index is less than two. This
is a natural condition in the context of jump filtering (see, e.g.,
\cite{Mancini2011} in the context of volatility estimation).
\end{remark}
%
\begin{remark}
To compare the finite to the infinite activity setting let us contrast
Assumption~\ref{ajumps}(ii) with Assumption \ref{aj}(ii). Both assumptions
control the behavior of small jumps. When the L\'{e}vy measure is finite
such that $\mu(\ddd x) = \lambda F(\ddd x)$ for some probability distribution
$F$ we can rewrite (\ref{smallj}) as
\[
\lambda\int_{-\dn^\beta}^{\dn^\beta} x^2 F (\ddd x) = \mathrm{O}
\bigl(\dn^\beta\bigr),
\]
since $\alpha= 0$ in this case. At the same time Assumption \ref{ajumps}(ii)
dominates Assumption \ref{aj}(ii) in the sense that
\[
\int_{-\dn^\beta}^{\dn^\beta} x^2 \mu(\ddd x) \leq\mu
\bigl(\bigl(-\dn^\beta, \dn^\beta\bigr)\bigr) = \lambda F\bigl(
\bigl(-\dn^\beta,\dn^\beta\bigr)\bigr).
\]
Hence, if $F((-\dn^\beta, \dn^\beta)) = \RMO (\dn^\beta)$ then
Assumption \ref {ajumps}(ii) implies Assumption \ref{aj}(ii) such that
a direct comparison becomes possible when $F$ has a bounded Lebesgue
density as in Example~\ref{densityex} below.
\end{remark}
%
\begin{lemma}\label{elemma}
If the L\'{e}vy measure $\mu_{|[-1,1]}$ of $M$ is symmetric around zero,
then Assumption~\ref{aj}\textup{(iii)} holds.
\end{lemma}
\begin{pf}
Assumption \ref{aj}(iii) is equivalent to the symmetry of the distribution
of $M$ restricted to $(-\delta,\delta)$ such that the statement
follows, since an infinitely divisible distribution is symmetric if and
only if its L\'evy--Khintchine triplet is of the form $(0,\sigma^2,\mu
)$ with $\mu$ being symmetric (cf.~\cite{Sato99}).
\end{pf}
%
\begin{remark}
Assumption \ref{aj}(iii) is\vspace*{1pt} a symmetry condition on the distribution
$P^{\di M}_{|(-\varepsilon,\varepsilon)}$ of the increments of $M$
restricted to $(-\varepsilon,\varepsilon)$ for every $0< \varepsilon\leq
\delta$. A sufficient condition for Assumption~\ref{aj}(iii) in terms of the
L\'{e}vy measure of $L$ was given in Lemma \ref{elemma}. But it is easy
to see that this is not a necessary condition. The main point here is
that $P^{\di M}_{|(-\varepsilon,\varepsilon)}$ is not infinitely divisible anymore.

Since our method does not depend on the choice of the maximal jump size
in the definition of $M$ in (\ref{jnotation}), it follows that the
condition given in Lemma \ref{elemma} can be relaxed to $\mu
_{|(-\varepsilon,\varepsilon)}$ being symmetric for some $\varepsilon>0$ by
redefining $M$ to have a L\'{e}vy measure supported on
$(-\varepsilon,\varepsilon)$.
\end{remark}

The main result of this chapter is the following central limit theorem
for the drift estimator with jump filter.
%
\begin{theorem} \label{discreteeff}
Suppose\vspace*{1pt} that Assumption \ref{aj} holds and $\sigma^2 > 0$. If there
exists $\beta\in(0,1/2)$ such that $T_n \Delta_n^{1 - 2\beta\wedge
({1/2})} = \RMo(1)$ as $n \to\infty$ then $v_n = \Delta_n^\beta$ yields
\[
T_n^{1/2} (\bar a_n - a) \stackrel{
\mathcal{D}} {\longrightarrow}N\bigl(0, \sigma^2 E_a
\bigl[X_\infty^2\bigr]^{-1}\bigr).
\]
The estimator is asymptotically efficient.
\end{theorem}
%
\begin{example}\label{densityex}
Let $L = W + J$, where $J$ is a compound Poisson process
\[
J_t = \sum_{i=1}^{N_t}
Y_i
\]
such that $Y_i \sim F$ are i.i.d. and $N_t$ is a Poisson process with
intensity $\lambda$. Suppose that $F$ has a bounded Lebesgue density
$f$. Then
\[
\int_{-v}^v x^2 \mu(\ddd x) = \lambda
\int_{-v}^v x^2 f(x) \dd x \leq C
v^3
\]
for $C >0$ such that for $L$ Assumption \ref{aj}(i) holds for every
$\alpha\in[0,2)$.

More generally every L\'evy process with Blumenthal--Getoor index less
than two fulfills Assumption \ref{aj}(i). This includes all L\'evy
processes commonly used in applications like (tempered) stable, normal
inverse Gaussian, variance gamma and also gamma processes.
\end{example}

\subsection{Proofs}
Asymptotic efficiency of $\bar a_n$ follows from the first statement of
Theorem \ref{discreteeff} together with Theorem \ref{consistencyThm}
such that it remains to prove the asymptotic normality result.
We will divide the proof of Theorem \ref{discreteeff} into several
lemmas. In the proofs in this section, constants may change from line
to line or even within one line without further notice.

\subsubsection{A moment bound}

In this section, we derive a moment bound for short time increments of
pure jump L\'{e}vy processes. Set
\[
f(x) = \cases{ x^2, &\quad if $|x| \leq1$,
\cr
0, &\quad if $|x| > 2$,}
\]
and $f(x) \in[0,2]$ for $|x| \in(1,2]$ such that $f \in C^\infty(\R
)$. We scale $f$ to be supported on $[-v,v]$ by
%
\begin{equation}
\label{fscaling} f^v (x) = v^2 f(x/v).
\end{equation}

\begin{proposition} \label{lemmoments}
Let $(M_t)_{t\geq0}$ be a pure jump L\'{e}vy process with L\'evy measure
$\mu$ such that $\operatorname{supp} (\mu) \subset[-1,1]$ and
Assumption \ref{aj}\textup{(i)} and \textup{(ii)} hold.
Then for all $\beta\in(0,\frac{1}{2})$ we obtain
\[
E \bigl[f^{t^\beta}(M_t) \bigr] = \mathrm{O} \bigl( t^{1 + \beta(2 - \alpha
)}
\bigr)
\]
as $t \downarrow0$.
\end{proposition}
%
\begin{remark}
The estimate in Proposition \ref{lemmoments} gives actually a bound
for the Markov generator of $M$ on the smooth test function $f^v$.
\end{remark}

\begin{pf*}{Proof of Proposition \ref{lemmoments}}
Let $P^{M_t}$ denote the distribution of $M_t$. We apply Plancherel's
identity to obtain
\[
E\bigl[f^{t^\beta}(M_t)\bigr] = \int_\R
f^{t^\beta}(x) P^{M_t} (\ddd x) = (2 \uppi )^{-1} \int
_\R\mathfrak{F} f^{t^\beta}(u) \overline{
\phi_t (u)} \dd u,
\]
where $\mathfrak{F} f = \int_{\R} \RMe^{\mathrm{i}ux} f(x) \dd x$ denotes the
Fourier transform of $f$ and the characteristic function of $M$ satisfies
\[
\phi_t (u) = \exp\biggl( t \int_{-1}^1
\bigl(\RMe^{\mathrm{i}ux} - 1 - \mathrm{i}ux\bigr) \mu(\ddd x) \biggr).
\]
Let us rewrite $\phi_t$ as the linearization of the exponential at
zero plus a remainder $R$:
\[
\phi_t (u) = 1 + \psi_t (u) + R(t,u)
\]
with
\[
\psi_t (u) = t \int_{-1}^1
\bigl(\RMe^{\mathrm{i}ux} - 1 - \mathrm{i}ux\bigr) \mu(\ddd x).
\]
Then,
%
\begin{eqnarray}
\label{fMdecomp} E\bigl[f^{t^\beta}(M_t)\bigr] &=& (2
\uppi )^{-1} \int_\R\mathfrak{F} f^{t^\beta
}(u)
\bigl(1 + \overline{\psi_t (u)} + \overline{R(t,u)}\bigr) \dd u
\nonumber\\[-8pt]\\[-8pt]
&=& (2 \uppi )^{-1} \int_\R\mathfrak{F}
f^{t^\beta}(u) \overline{\psi_t (u)} \dd u + (2
\uppi )^{-1} \int_\R\mathfrak{F} f^{t^\beta}(u)
\overline{R(t,u)} \dd u.\nonumber
\end{eqnarray}
For the first term on the right-hand side, we obtain
%
\begin{eqnarray}
\label{firstterm} (2 \uppi )^{-1} \int_\R
\mathfrak{F} f^{t^\beta}(u) \overline{\psi_t (u)} \dd u &=& (2
\uppi )^{-1} t \int_{-1}^1 \int
_\R\mathfrak{F} f^{t^\beta}(u) \bigl(\RMe^{-\mathrm{i}ux}
-1 + \mathrm{i}ux\bigr) \dd u \mu(\ddd x)
\nonumber
\\
&=& t \int_{-1}^1 \biggl(f^{t^\beta}(x) +
(2 \uppi )^{-1} \int_\R\mathfrak{F} \bigl(
\bigl(f^{t^\beta}\bigr)' \bigr) (u) x \dd u \biggr) \mu(\ddd x)\quad
\\
&=& t \int_{-1}^1 f^{t^\beta}(x) \mu(\ddd x) =
t \mathrm{O}\bigl(t^{\beta(2-\alpha)}\bigr)\nonumber
\end{eqnarray}
by Assumption \ref{aj}(i) and since
\[
\int_\R\mathfrak{F} \bigl( \bigl(f^{t^\beta}
\bigr)' \bigr) (u) \dd u = \bigl( f^{t^\beta}
\bigr)' (0) = 0.
\]
It remains to bound the second addend in (\ref{fMdecomp}). For $\operatorname{Re}(z)
\leq0$ observe that
%
\begin{equation}
\label{expest} \biggl\llvert\frac{\RMe^z - z -1}{z^2} \biggr\rrvert\leq C
\end{equation}
for constant $C > 0$, since for $|z| \geq1$
\[
\frac{|\RMe^z - z -1|}{z^2} \leq2 + \frac{|z|}{z^2} \leq3.
\]
Whereas on the half disk $\{|z| <1, \operatorname{Re}(z) \leq0 \}$ the continuous
function $|\RMe^z -z -1|$ is bounded and $z^2$ is bounded except for the
singularity at the origin, but at zero we know that $|\RMe^z -z -1| =
\RMO(z^2)$, that is,
\[
\frac{|\RMe^z - z -1|}{z^2} \leq C < \infty
\]
on $\{|z| <1, \operatorname{Re}(z) \leq0 \}$. Theorem 1.2.5 in \cite{Kappus2012}
implies that $|\psi_t (u)| \leq C t |u|^{\alpha}$ such that
\[
\bigl|R(t,u)\bigr|= \bigl|\RMe^{\psi_t (u)} - \psi_t (u) -1\bigr| \leq\bigl\llvert
\psi_t (u)\bigr\rrvert^2 \leq C t^2
|u|^{2 \alpha},
\]
where we used (\ref{expest}) and that for every characteristic
function $|{\exp(\psi_t (u))}| = \phi_t (u) \leq1$ holds. Hence, we obtain
%
\begin{equation}
\label{Restimate} \biggl\llvert\int_\R\mathfrak{F}
\bigl(f^{t^\beta} \bigr) (u) \overline{R(t,u)} \dd u \biggr\rrvert
\leq C
t^2 \int_\R\bigl\llvert\mathfrak{F}
\bigl(f^{t^\beta} \bigr) (u) \bigr\rrvert|u|^{2 \alpha} \dd u.
\end{equation}
Therefore, it remains to bound $\int_\R\llvert \mathfrak{F}
f^{t^\beta}(u) \rrvert |u|^{2 \alpha} \dd u$ in $t$. From (\ref
{fscaling}) and the scaling property of the Fourier transform it
follows that
\[
\mathfrak F \bigl(f^v\bigr) (u) = v^3 \mathfrak F
\bigl(v^{-1} f(x/v) \bigr) (u) = v^3 \mathfrak F (f) (vu).
\]
Since $f \in C^\infty(\R)$, we obtain $|\mathfrak F (f) (u)| \leq C_m
|u|^{-m}$ such that
\[
\bigl| \mathfrak F \bigl(f^v\bigr) (u)\bigr| \leq C_m
v^{3-m} |u|^{-m}
\]
for all $u \in\R$ and $m,v >0$. Then
\[
h(v,u) = \bigl|\mathfrak F \bigl(f^v\bigr) (u)\bigr| |u|^{2\alpha} \leq
C_m v^{3-m} |u|^{2\alpha-m}.
\]
If
%
\begin{equation}
\label{alphamrel} 2\alpha+1 < m
\end{equation}
holds then $h(v,\cdot) \in L^1 (\R)$ for all $v \in(0,1)$. Setting
$v = t^\beta$ yields
\[
t^2 \int_\R\bigl\llvert\mathfrak{F}
\bigl(f^{t^\beta} \bigr) (u) \bigr\rrvert|u|^{2 \alpha} \dd u \leq
C_m t^{(3-m) \beta+2}
\]
for all $m >0$. Since the first term in (\ref{fMdecomp}) is of the
order $\RMO (t^{1+ \beta(2-\alpha)} )$, we choose $m$ such that
\[
(3-m) \beta+ 2 \geq1 + \beta(2-\alpha) \quad\Leftrightarrow\quad m \leq1+
\beta^{-1} + \alpha.
\]
Together with (\ref{alphamrel}) this leads to the condition
\[
2 \alpha+1 < 1+ \beta^{-1} + \alpha\quad\Leftrightarrow\quad\alpha<
\beta^{-1},
\]
which due to $\alpha\in(0,2)$ always holds for $\beta\in(0,1/2)$.
Hence, we obtain
\[
\biggl\llvert\int_\R\mathfrak{F} \bigl(f^{t^\beta}
\bigr) (u) \overline{R(t,u)} \dd u \biggr\rrvert= \mathrm{O} \bigl( t^{1+ \beta
(2-\alpha)}
\bigr).
\]
Together with (\ref{fMdecomp}) and (\ref{firstterm}) this yields finally
\[
E\bigl[f^{t^\beta}(M_t)\bigr] = t \mathrm{O}\bigl(t^{\beta(2-\alpha)}
\bigr).
\]
\upqed\end{pf*}

\subsubsection{Approximating the continuous martingale part}
The main step is to show that the continuous martingale part can be
approximated by summing only the increments that are below the
threshold $v_n$. We will use throughout the notation from~(\ref{jnotation}).
%
\begin{lemma} \label{jumpfilter}
Suppose that the assumptions of Theorem \ref{discreteeff} hold, then
\[
T_n^{-1/2} \sum_{i=0}^{n-1}X_{t_i}
\bigl(\di X \1_{\{|\di X| \leq v_n\}} - \di X^c\bigr) \stackrel{p} {
\longrightarrow}0 \qquad\mbox{as $n \to\infty$}.
\]
\end{lemma}

\begin{pf}
Let us consider the following decomposition where $\tilde X = X_0 + W +
D + U$
\begin{eqnarray*}
&&
T_n^{-1/2} \sum_{i=0}^{n-1}X_{t_i}
\bigl(\di X \1_{\{|\di X| \leq v_n\}} - \di X^c\bigr)
\\
&&\quad= T_n^{-1/2} \sum_{i=0}^{n-1}X_{t_i}
\bigl(\di\tilde X \1_{\{|\di X| \leq
v_n\}} - \di X^c\bigr)
+T_n^{-1/2} \sum_{i=0}^{n-1}X_{t_i}
\di M \1_{\{|\di X| \leq v_n\}}
\\
&&\quad= T_n^{-1/2} \sum_{i=0}^{n-1}X_{t_i}
\bigl(\di\tilde X \1_{\{|\di\tilde
X| \leq
2v_n\}} - \di X^c\bigr)
\\
&&\qquad{} + T_n^{-1/2} \sum_{i=0}^{n-1}X_{t_i}
\di\tilde X (\1_{\{|\di
X| \leq v_n\}} - \1_{\{|\di\tilde X| \leq2v_n\}} ) +T_n^{-1/2}
\sum_{i=0}^{n-1}X_{t_i} \di M
\1_{\{|\di X| \leq v_n\}}
\\
&&\quad= S_1^n + S_2^n +
S_3^n.
\end{eqnarray*}
Observe that the term $S_1^n$ already appeared in Lemma \ref{JApp} and
$\tilde X = X_t - M_t$ is a process with finite jump activity. A
careful analysis of the proof of Lemma \ref{JApp} reveals that the
same estimates apply to $S_1^n$ such that we conclude that $S_1^n$
converges to zero in probability when $n \to\infty$. Let us prove
next convergence of
\[
S_2^n = T_n^{-1/2} \sum
_{i=0}^{n-1}X_{t_i} \di\tilde X (-
\1_{\{
|\di X| >
v_n, |\di\tilde X| \leq2v_n\}} + \1_{\{|\di X| \leq v_n, |\di\tilde
X| > 2v_n\}} ).
\]
Let us prove next that the contribution of the second indicator
function on the right-hand side tends to zero in probability:
%
\begin{eqnarray}
\label{s2est}
&& P \Biggl(T_n^{-1/2} \sum
_{i=0}^{n-1}|X_{t_i} \di\tilde X|
\1_{\{|\di
X| \leq
v_n, |\di\tilde X| > 2v_n\}} > 0 \Biggr)
\nonumber\\
&&\quad= P \Biggl( \bigcup_{i=0}^{n-1} \bigl\{|\di X|
\leq v_n, |\di\tilde X| > 2v_n\bigr\} \Biggr) \\
&&\quad\leq\sum
_{i=0}^{n-1}P \bigl(|\di X| \leq v_n,
|\di\tilde X| > 2v_n\bigr).\nonumber
\end{eqnarray}
When $|\di\tilde X|> 2 v_n$ then with high probability $|\di U| > 0$,
since by Lemma \ref{pbound} we obtain
%
\begin{eqnarray}
\label{xuestimate} \sum_{i=0}^{n-1}P\bigl(|\di
\tilde X|> 2 v_n, |\di U| =0\bigr) &\leq&\sum_{i=0}^{n-1}P\bigl(|
\di W + \di D| > 2 v_n \bigr) \nonumber\\[-8pt]\\[-8pt]
&=& \mathrm{O}\bigl(T_n \dn^{1-2\beta}
\bigr).\nonumber
\end{eqnarray}
This together with (\ref{s2est}) and that fact that on $\{|\di X| \leq
v_n, |\di\tilde X| > 2v_n\}$ necessarily $|\di M| > v_n$ implies that
\begin{eqnarray*}
&&
P \Biggl(T_n^{-1/2} \sum_{i=0}^{n-1}|X_{t_i}
\di\tilde X| \1_{\{|\di
X| \leq
v_n, |\di\tilde X| > 2v_n\}} > 0 \Biggr)
\\
&&\quad\leq\sum_{i=0}^{n-1}P\bigl(|\di U| \neq0\bigr) P\bigl(|
\di M| > v_n\bigr) + \mathrm{O} \bigl(T_n \Delta_n^{1 - 2\beta
}
\bigr)
\\
&&\quad= \mathrm{O} \bigl(T_n \Delta_n v_n^{-2}
\bigr) + \mathrm{O} \bigl(T_n \Delta_n^{1 - 2\beta}\bigr) = \mathrm{O}
\bigl(T_n \Delta_n^{1 - 2\beta}\bigr),
\end{eqnarray*}
where we used Markov's inequality and independence of $U$ and $M$. The
remaining term in $S_2^n$ is
\[
T_n^{-1/2} \sum_{i=0}^{n-1}X_{t_i}
\di\tilde X \1_{\{|\di X| > v_n,
|\di\tilde
X| \leq2v_n\}}.
\]
Let us prove that on $\{ |\di\tilde X| \leq2 v_n\}$ the contribution
of $U$ is negligible:
%
\begin{eqnarray}
\label{unegligible}
&&
T_n^{-1/2} \sum
_{i=0}^{n-1}P\bigl(|\di\tilde X| \leq2 v_n, |
\di U| > 0\bigr)
\nonumber
\\
&&\quad= T_n^{-1/2} \sum_{i=0}^{n-1}
\bigl( P\bigl(|\di\tilde X| \leq2 v_n, \di N = 1\bigr) + \mathrm{O}\bigl(
\dn^2\bigr) \bigr)
\\
&&\quad\leq T_n^{-1/2} \sum_{i=0}^{n-1}
\bigl( P\bigl(|\di W + \di D| > 1 - 2 v_n\bigr) + \mathrm{O}\bigl(\dn^
2
\bigr) \bigr) = \mathrm{O}\bigl(T_n^{1/2} \dn\bigr)
\nonumber
\end{eqnarray}
as $n \to\infty$, where $N$ denotes the counting process that counts
the jumps of $U$ and the last step follows from Lemma \ref{pbound}.
Hence, we can assume that $\di U =0$ on $\{ |\di\tilde X| \leq2 v_n\}
$ and so $\di\tilde X = \di W + \di D$.
For $T_n \Delta_n^{1/2-\beta/2} = \RMo(1)$ it follows from Lemma \ref
{lem1} that as $n \to\infty$
\[
T_n^{-1/2} \sum_{i=0}^{n-1}X_{t_i}
\di D \1_{\{|\di X| > v_n, |\di
\tilde X| \leq
2v_n\}} \stackrel{p} {\longrightarrow}0.
\]
We have decomposed $S_2^n$ into a term that converges to $0$ in
probability and a remainder:
\[
S_2^n= T_n^{-1/2} \sum
_{i=0}^{n-1}X_{t_i} \di W \1_{\{|\di X| > v_n,
|\di\tilde
X| \leq2v_n\}}
+ \RMo_p (1).
\]
For the remainder let us observe that by Lemma \ref{lemj}, we obtain
\begin{eqnarray*}
S_2^n &=& T_n^{-1/2} \sum
_{i=0}^{n-1}X_{t_i} \di W \1_{\{|\di X| > v_n,
|\di
\tilde X| \leq2v_n\}}
+ \RMo_p(1)
\\
&=& T_n^{-1/2} \sum_{i=0}^{n-1}X_{t_i}
\di W \1_{\{|\di W + \di D + \di
M| > v_n,
|\di W + \di D| \leq2v_n\}} + \RMo_p(1)
\\
&=& T_n^{-1/2} \sum_{i=0}^{n-1}X_{t_i}
\di W \1_{\{|\di M| > v_n\}} + \RMo_p(1)
\end{eqnarray*}
Markov's inequality yields $P(|\di M| > v_n) \leq\Delta_n^{1/2 -
\beta}$. Independence of $X_{t_i}$, $\di W$ and $\di M$ leads to
\[
E \Biggl[ T_n^{-1/2} \sum_{i=0}^{n-1}X_{t_i}
\di W \1_{\{|\di M| >
v_n\}} \Biggr] =0
\]
and
\begin{eqnarray*}
&&
E \Biggl[ \Biggl(T_n^{-1/2} \sum
_{i=0}^{n-1}X_{t_i} \di W \1_{\{|\di
M| > v_n\}}
\Biggr)^2 \Biggr] \\
&&\quad\leq T_n^{-1} E \Biggl[\sum
_{i=0}^{n-1}X_{t_i}^2 (
\di W)^2 \1 _{\{|\di M| > v_n\}} \Biggr]
\\
&&\qquad{}+ T_n^{-1} E \biggl[ \sum
_{ i \neq j} X_{t_i} \di W \1_{\{|\di M| >
v_n\}}
X_{t_j} \Delta_j W \1_{\{|\Delta_j M| > v_n\}} \biggr].
\end{eqnarray*}
Since $X_{t_i}, X_{t_j}, \di W, \di M, \Delta_j M\perp\Delta_j W$
for $i < j$ the off-diagonal elements are centered,
\[
E \biggl[ \sum_{ i \neq j} X_{t_i} \di W
\1_{\{|\di M| > v_n\}} X_{t_j} \Delta_j W \1_{\{|\Delta_j M| > v_n\}}
\biggr] =0
\]
and the diagonal elements can be estimated by
\begin{eqnarray*}
T_n^{-1} E \Biggl[\sum_{i=0}^{n-1}X_{t_i}^2
(\di W)^2 \1_{\{|\di M| >
v_n\}
} \Biggr] &\leq& T_n^{-1}
\Delta_n \sum_{i=0}^{n-1}E
\bigl[X_{t_i}^2\bigr] P\bigl(|\di M| > v_n\bigr)
\\
&\leq& \sup_i \bigl\{ E\bigl[X_{t_i}^2
\bigr] \bigr\} \Delta_n^{1/2-\beta} \to0
\end{eqnarray*}
as $n \to\infty$. The last step is to show that $S_3^n$ tends to zero
in probability as $n \to\infty$. As in (\ref{unegligible}) it
follows that on $|\di X| \leq v_n$ we can assume that $\di U = 0$. Now
\begin{eqnarray*}
\sum_{i=0}^{n-1}P\bigl(|\di X| \leq
v_n, \di U =0\bigr) &\leq& \sum_{i=0}^{n-1}P\bigl(|
\di W + \di D +\di M| \leq v_n, |\di M| \leq2v_n\bigr)
\\
&&{}+ \sum_{i=0}^{n-1}P\bigl(|\di W + \di D +\di
M| \leq v_n, |\di M| > 2v_n\bigr).
\end{eqnarray*}
The second addend vanishes, since by Lemma \ref{pbound} we obtain
\begingroup
\abovedisplayskip=7pt
\belowdisplayskip=7pt
\begin{eqnarray*}
&&
\sum_{i=0}^{n-1}P\bigl(|\di W + \di D +\di M|
\leq v_n, |\di M| > 2v_n\bigr)
\\[-2pt]
&&\quad\leq\sum_{i=0}^{n-1} P\bigl(|\di W + \di D | >
v_n\bigr) = \mathrm{O}\bigl(T_n \dn^{1-2 \beta}\bigr).
\end{eqnarray*}
Thus, $S_3^n$ can be rewritten as
%
\begin{equation}
\label{S3decomp} S_3^n = T_n^{-1/2}
\sum_{i=0}^{n-1}X_{t_i} \di M
\1_{\{|\di X| \leq
v_n, |\di M|
\leq2 v_n\}} + \RMo_p(1).
\end{equation}
The convergence of the remaining term in $S_3^n$ is dominated by the
behavior of $\di M$ around the threshold, that is, we prove next that
\begin{eqnarray*}
&&
T_n^{-1/2} \sum_{i=0}^{n-1}X_{t_i}
\di M \1_{\{|\di X| \leq v_n, |\di
M| \leq2
v_n\}}\\[-2pt]
&&\quad = T_n^{-1/2} \sum
_{i=0}^{n-1}X_{t_i} \di M \1_{\{ |\di M| \leq
2 v_n\}}
+ \RMo_p(1).
\end{eqnarray*}
Indeed,
\begin{eqnarray*}
&&
T_n^{-1/2} \sum_{i=0}^{n-1}X_{t_i}
\di M ( \1_{\{ |\di M| \leq2 v_n\}
} - \1_{\{
|\di X| \leq v_n, |\di M| \leq2 v_n\}})
\\[-2pt]
&&\quad= T_n^{-1/2} \sum_{i=0}^{n-1}X_{t_i}
\di M \1_{\{|\di X| > v_n, |\di
M| \leq2
v_n\}}.
\end{eqnarray*}
That last term tends to zero in probability will be shown in the proof
of Lemma \ref{lemj} below following equation (\ref{xmcon}). To finish
the proof, we demonstrate that the first addend on the right-hand side
of (\ref{S3decomp}) vanishes asymptotically.
Since $X_{t_i}, X_{t_j}, \di M \perp\Delta_j M$ for $i < j$ the
off-diagonal elements vanish by Assumption \ref{aj}(i) such that
\[
E \biggl[ \sum_{ i \neq j} X_{t_i} \di M
\1_{\{ |\di M| \leq2 v_n\}} X_{t_j} \Delta_j M\1_{\{ |\Delta_j M| \leq2
v_n\}}
\biggr] =0
\]
and the diagonal elements can by Proposition \ref{lemmoments} be
estimated by
\begin{eqnarray*}
T_n^{-1} E \Biggl[\sum_{i=0}^{n-1}X_{t_i}^2
(\di M)^2 \1_{\{ |\di M|
\leq2 v_n\}
} \Biggr] &\leq& T_n^{-1}
\sum_{i=0}^{n-1}E\bigl[X_{t_i}^2
\bigr] E\bigl[ (\di M)^2 \1 _{\{ |\di
M| \leq2 v_n\}}\bigr]
\\[-2pt]
&\leq& \sup_i \bigl\{E\bigl[X_{t_i}^2
\bigr]\bigr\} v_n^\alpha\to0
\end{eqnarray*}
as $n \to\infty$.\vadjust{\goodbreak}
\end{pf}
\endgroup

\subsubsection{Approximation of the drift}
The next step is to show that the drift component of $X$ is in the
limit not affected by the cut-off.
%
\begin{lemma}\label{lem1}
If the assumptions of Theorem \ref{discreteeff} are fulfilled then
\[
T_n^{-1/2} \sum_{i=0}^{n-1}X_{t_i}
(\di D \1_{\{|\di X| \leq v_n\}} - \di D) \stackrel{p} {\longrightarrow
}0 \qquad\mbox{as $n
\to\infty$}.
\]
\end{lemma}
\begin{pf}
We rewrite the sum as follows:
\[
T_n^{-1/2} \sum_{i=0}^{n-1}X_{t_i}
(\di D \1_{\{|\di X| \leq v_n\}} - \di D) = T_n^{-1/2} \sum
_{i=0}^{n-1}X_{t_i} \di D
\1_{\{|\di X| > v_n\}}.
\]
Next, we decompose $\di D$ as follows
\[
\di D = -a \biggl( \int_{t_i}^{t_{i+1}}
(X_s - X_{t_i}) \dd s + \di X_{t_i} \biggr)
\]
such that by Lemma \ref{lemj} below
%
\begin{eqnarray}
\label{sum1} T_n^{-1/2}\sum_{i=0}^{n-1}X_{t_i}
\di D \1_{\{|\di X| > v_n\}} &=& -a T_n^{-1/2}\sum
_{i=0}^{n-1}X_{t_i} \int
_{t_i}^{t_{i+1}} (X_s - X_{t_i})
\dd s \1 _{\{|\di J| > v_n\}}
\nonumber\\[-8pt]\\[-8pt]
&&{}-a T_n^{-1/2} \sum_{i=0}^{n-1}X_{t_i}^2
\di\1_{\{|\di J| > v_n\}} + \RMo_p (1).\nonumber
\end{eqnarray}
For the second term, we obtain by Markov's inequality and from $v_n =
\Delta_n^\beta$ that
\[
E \Biggl[ \Biggl\llvert\sum_{i=0}^{n-1}X_{t_i}^2
\di\1_{\{|\di J| > v_n\}
}\Biggr\rrvert\Biggr] \leq\sum_{i=0}^{n-1}
\di E\bigl[X_{t_i}^2\bigr] P \bigl(|\di J| > v_n\bigr)
\leq C T_n \Delta_n^{1-2\beta}
\]
and so for $T_n^{1/2} \Delta_n^{1-2 \beta} = \RMo(1)$ it follows that
\[
T_n^{-1/2} \sum_{i=0}^{n-1}X_{t_i}
\di\1_{\{|\di J| > v_n\}} = \RMo_p (1).
\]
For the first sum on the right-hand side of (\ref{sum1}), we obtain by
H\"older's inequality and independence of $X_{t_i}$ and $\di J$
\begin{eqnarray*}
&&
E \biggl[ \biggl\llvert X_{t_i} \int_{t_i}^{t_{i+1}}
(X_s - X_{t_i}) \dd s \1 _{\{|\di J| > v_n\}} \biggr\rrvert
\biggr]
\\
&&\quad\leq E \biggl[ \biggl( \int_{t_i}^{t_{i+1}}
(X_s - X_{t_i}) \dd s \biggr)^2
\biggr]^{1/2} P\bigl(|\di J| > v_n\bigr)^{1/2} E
\bigl[X_{t_i}^2\bigr]^{1/2} \\
&&\quad= \mathrm{O}\bigl(\Delta
_n^{3/2}v_n^{-1/2}\bigr)
\end{eqnarray*}
such that for $T_n^{1/2} \Delta_n^{1/2 - \beta/2} = \RMo(1)$ we can
conclude that
\[
T_n^{-1/2} \sum_{i=0}^{n-1}X_{t_i}
\int_{t_i}^{t_{i+1}} (X_s -
X_{t_i}) \dd s \1 _{\{|\di J| > v_n\}} = \RMo_p(1).
\]
\upqed\end{pf}

\subsubsection{Identifying the jumps}

In the following, we will show that the increments of $X$ that are
larger than the threshold $v_n$ are dominated by the jump component.
%
\begin{lemma} \label{lemj}
\[
T_n^{-1/2} \sum_{i=0}^{n-1}X_{t_i}
\di X (\1_{\{|\di X| \leq v_n\}} - \1_{\{|\di
M| \leq2 v_n\}}) \stackrel{p} {\longrightarrow}0
\qquad\mbox{as $n \to\infty$}.
\]
\end{lemma}
\begin{pf}
Observe that
%
\begin{eqnarray}
\label{decomp1}
&&
T_n^{-1/2} \sum
_{i=0}^{n-1}X_{t_i} \di X (
\1_{\{|\di X| \leq v_n\}} - \1_{\{|\di
J| \leq2v_n\}})
\\
&&\quad= T_n^{-1/2} \sum_{i=0}^{n-1}X_{t_i}
\di X \1_{\{
|\di X| \leq
v_n, |\di J| > 2 v_n\}}
\\
&&\qquad{}- T_n^{-1/2} \sum_{i=0}^{n-1}X_{t_i}
\di X \1_{\{|\di X| > v_n, |\di
J| \leq2
v_n\}}.
\nonumber
\end{eqnarray}
We shall prove in Lemma \ref{lemJ} below that
%
\begin{equation}
\label{con1} T_n^{-1/2} \sum_{i=0}^{n-1}X_{t_i}
\di X \1_{\{|\di X| \leq v_n, |\di
J| > 2 v_n\}
} \stackrel{p} {\longrightarrow}0.
\end{equation}
In the next step, we show that the contribution of $U$ is negligible,
since by independence of $\di W$, $\di M$, $\di U$ and $X_{t_i}$ it
follows that
%
\begin{eqnarray}
\label{decomp2}
E \Biggl[ \Biggl\llvert\sum_{i=0}^{n-1}X_{t_i}
\di X \1_{\{|\di J| \leq
2v_n, |\di U|
\neq0\}}\Biggr\rrvert\Biggr] &\leq&\sum
_{i=0}^{n-1}E\bigl[|X_{t_i}|\bigr] E\bigl[|\di W|\bigr] P\bigl( |\di
U| \neq0\bigr)
\nonumber\\
&&{}+ \sum_{i=0}^{n-1}E\bigl[|X_{t_i}|\bigr]
E\bigl[|\di J| \1_{\{|\di J| \leq2v_n,
|\di U| \neq
0\}}\bigr]\qquad\quad  \\
&&{}+ \sum_{i=0}^{n-1}E\bigl[|X_{t_i}
\di D| \1_{\{ |\di U| \neq0\}}\bigr].\nonumber\quad
\end{eqnarray}
Now $U$ is a compound Poisson process with intensity $\mu( \R
\setminus[-1,1])< \infty$ such that $P(\di U \neq0) = \RMO(\dn)$. We
obtain for first addend on the right-hand side
\[
T_n^{-1/2} \sum_{i=0}^{n-1}E
\bigl[|\di W|\bigr] E\bigl[|X_{t_i}|\bigr] P(\di U \neq0) = \mathrm{O}\bigl( T_n^{1/2}
\dn^{1/2}\bigr).
\]
We split the second term into the contribution by $U$ and $M$ such that
\begin{eqnarray*}
&&
T_n^{-1/2} \sum_{i=0}^{n-1}E\bigl[|X_{t_i}|\bigr]
E\bigl[|\di J|\1_{\{|\di J| \leq
2v_n,\di U
\neq0\}}\bigr]
\\
&&\quad= T_n^{-1/2} \sum_{i=0}^{n-1}E\bigl[|X_{t_i}|\bigr]
E\bigl[|\di M|\bigr] E[\1_{\{|\di J|
\leq
2v_n,\di U \neq0\}}]
\\
&&\qquad{}+T_n^{-1/2} \sum_{i=0}^{n-1}E\bigl[|X_{t_i}|\bigr]
E\bigl[|\di U|\1_{\{
|\di J| \leq2v_n, \di U \neq0\}}\bigr].
\end{eqnarray*}
The first sum is of order
\[
T_n^{-1/2} \sum_{i=0}^{n-1}E\bigl[|X_{t_i}|\bigr]
E\bigl[|\di M|\bigr] E[\1_{\{\di U \neq
0\}}] = \mathrm{O}\bigl(T_n^{1/2}
\dn^{1/2}\bigr).
\]
H\"oldern's inequality and independence of $M$ and $U$ lead to the
following estimate for the second sum:
\[
T_n^{-1/2} E \sum_{i=0}^{n-1}|X_{t_i}
\di U\1_{\{|\di J| \leq2v_n,\di
U \neq0\}}| = \mathrm{O}\bigl(T_n^{1/2}
\dn^{1/2}\bigr).
\]
To prove convergence of the last addend in (\ref{decomp2}), we rewrite
$\di D$ as follows
%
\begin{equation}
\label{Ddecomposition} \di D = -a \biggl( \int_{t_i}^{t_{i+1}}
(X_s - X_{t_i}) \dd s + \di X_{t_i} \biggr)
\end{equation}
and so
\begin{eqnarray*}
T_n^{-1/2} E\Biggl[\Biggl| \sum_{i=0}^{n-1}X_{t_i}
\di D\1_{\{\di U \neq0\}}\Biggr|\Biggr] &\leq& T_n^{-1/2}a E\Biggl[\Biggl|
\sum_{i=0}^{n-1}X_{t_i} \int
_{t_i}^{t_{i+1}} (X_s - X_{t_i})
\dd s \1_{\{\di U \neq0\}}\Biggr|\Biggr]
\\
&&{}+ T_n^{-1/2} E\Biggl[\Biggl| \sum_{i=0}^{n-1}X_{t_i}^2
\dn\1_{\{\di U \neq0\}}\Biggr|\Biggr].
\end{eqnarray*}
The first term on the right-hand side gives by using H\"older's inequality
\begin{eqnarray*}
&&
T_n^{-1/2}a E \Biggl[\Biggl\llvert\sum
_{i=0}^{n-1}X_{t_i} \int
_{t_i}^{t_{i+1}} (X_s - X_{t_i})
\dd s \1_{\{\di U \neq0\}} \Biggr\rrvert\Biggr]
\\
&&\quad\leq T_n^{-1/2}a \sum_{i=0}^{n-1}E
\bigl[X_{t_i}^2\bigr]^{1/2} E \biggl[ \biggl(\int
_{t_i}^{t_{i+1}} (X_s - X_{t_i})
\dd s \biggr)^2 \biggr]^{1/2} \1_{\{
\di U \neq0\}}
\\
&&\quad=\mathrm{O}\bigl(T_n^{1/2} \dn\bigr).
\end{eqnarray*}
Hence, we obtain
\[
T_n^{-1/2} \sum_{i=0}^{n-1}X_{t_i}
\di D\1_{\{\di U \neq0\}} = \mathrm{O}_p \bigl(T_n^{1/2}
\dn^{1/2}\bigr)
\]
such that it follows that
\[
\sum_{i=0}^{n-1}X_{t_i} \di X
\1_{\{|\di J| \leq2v_n, |\di U| \neq0\}
} = \RMo_p(1).
\]
Since the contribution of $U$ is negligible, we obtain from (\ref
{decomp1}) and (\ref{con1}) that
\begin{eqnarray*}
&&T_n^{-1/2} \sum_{i=0}^{n-1}X_{t_i}
\di X (\1_{\{|\di X| \leq v_n\}} - \1_{\{|\di
M| \leq v_n\}})
\\
&&\quad= T_n^{-1/2} \sum_{i=0}^{n-1}X_{t_i}
\di X \1_{\{
|\di X| >
v_n, |\di J| \leq2 v_n, |\di U| = 0\}} + \RMo_p(1).
\end{eqnarray*}
Hence, it remains to prove
%
\begin{equation}
\label{xmcon} T_n^{-1/2} \sum_{i=0}^{n-1}X_{t_i}
\di X \1_{\{\di X| > v_n|\di M|
\leq2v_n,
|\di U| = 0\}} \stackrel{p} {\longrightarrow}0 \qquad\mbox{as $n \to
\infty$}.
\end{equation}
Observe that
\begin{eqnarray*}
&&\bigl\{ |\di M| \leq2 v_n, \di U =0, |\di X | > v_n \bigr\}
\\
&&\quad\subset\bigl\{ |\di M| \leq2 v_n, \di U = 0, |\di W + \di D| + |\di M|
> v_n\bigr\}
\\
&&\quad\subset\bigl\{ |\di W + \di D| > v_n/2 \bigr\} \cup\bigl\{ |\di M| \leq2
v_n, |\di M| > v_n/2 \bigr\}.
\end{eqnarray*}
Therefore, the last two steps will be to show that:
\begin{longlist}[(ii)]
\item[(i)] $
T_n^{-1/2} \sum_{i=0}^{n-1}X_{t_i} \di X \1_{\{|\di W + \di D| > v_n/2
\}} \stackrel{p}{\longrightarrow}0$,
\item[(ii)] $
T_n^{-1/2} \sum_{i=0}^{n-1}X_{t_i} \di X \1_{\{|\di M| > v_n/2, |\di
M| \leq2
v_n\}} \stackrel{p}{\longrightarrow}0$.
\end{longlist}
For the proof of these two convergences, we refer to Lemmas \ref{lemM}
and \ref{lemW}.\vadjust{\goodbreak}
\end{pf}

\begin{lemma} \label{lemJ}
\[
T_n^{-1/2} \sum_{i=0}^{n-1}X_{t_i}
\di X \1_{\{|\di X| \leq v_n, |\di
J| > 2 v_n\}
} \stackrel{p} {\longrightarrow}0 \qquad\mbox{as $n
\to\infty$}.
\]
\end{lemma}
\begin{pf}
On $\{ |\di X| \leq v_n, |\di J| \geq2 v_n\}$ we have
\[
\bigl||\di W + \di D | - |\di J|\bigr| \leq|\di X| \leq v_n.
\]
Hence, we necessarily have $|\di W + \di D| > v_n$, that is,
%
\begin{equation}
\label{sub1} \bigl\{ |\di X | \leq v_n, | \di J | > 2 v_n
\bigr\} \subset\bigl\{ |\di W + \di D | > v_n \bigr\}
\end{equation}
such that
%
\begin{equation}
\label{moments} P\bigl( |\di X| \leq v_n, | \di J | > 2 v_n
\bigr) \leq P \bigl( | \di W + \di D| > v_n \bigr)
= \mathrm{O}\bigl( \dn^{2- \beta}
\bigr).
\end{equation}
It follows from (\ref{sub1}) that
\begin{eqnarray*}
&&
T_n^{-1/2} \Biggl\llvert\sum
_{i=0}^{n-1}X_{t_i} \di X \1_{\{|\di X| \leq
v_n, |\di J|
> 2 v_n\}}
\Biggr\rrvert
\\[-2pt]
&&\quad\leq T_n^{-1/2} \sum_{i=0}^{n-1}|X_{t_i}
\di X | \1_{\{|\di W + \di
D|> v_n\}}
\\[-2pt]
&&\quad\leq T_n^{-1/2} \sum_{i=0}^{n-1}|X_{t_i}
\di W| \1_{\{|\di W + \di
D|> v_n\}} + T_n^{-1/2} \sum
_{i=0}^{n-1}|X_{t_i} \di D|
\1_{\{|\di W + \di D|> v_n\}
}
\\[-2pt]
&&\qquad{} + T_n^{-1/2} \sum_{i=0}^{n-1}|X_{t_i}
\di M| \1_{\{|\di W + \di
D|> v_n\}} + T_n^{-1/2} \sum
_{i=0}^{n-1}|X_{t_i} \di U|
\1_{\{|\di W + \di D|>
v_n\}}
\\[-2pt]
&&\quad= A_n^1+ \cdots+ A_n^4.
\end{eqnarray*}
For $A_n^1$ we find by (\ref{moments}), H\"older's inequality and
independence of $X_{t_i}$ and $\di W$ that
\begin{eqnarray*}
E\bigl[\bigl| A_n^1\bigr|\bigr] &\leq& T_n^{-1/2}
\dn^{1/2} \sum_{i=0}^{n-1}E
\bigl[X_{t_i}^2\bigr]^{1/2} P\bigl(|\di W + \di D|>
v_n\bigr)^{1/2}
\\[-2pt]
&=& \mathrm{O}\bigl(T_n^{1/2} \dn^{1/2 - \beta/2}\bigr).
\end{eqnarray*}
Using (\ref{Ddecomposition}), we obtain for $A_n^2$ that
\begin{eqnarray*}
E\bigl[\bigl| A_n^2\bigr|\bigr] &\leq& T_n^{-1/2}
a \sum_{i=0}^{n-1}E \biggl[ \biggl\llvert
X_{t_i} \int_{t_i}^{t_{i+1}} (X_s
- X_{t_i}) \dd s \biggr\rrvert\1_{\{|\di W + \di D|>
v_n\}} \biggr]
\\[-2pt]
&&{}+ T_n^{-1/2} \dn a \sum_{i=0}^{n-1}E
\bigl[ X_{t_i}^2 \1_{\{|\di W + \di
D|> v_n\}}\bigr]
\end{eqnarray*}
H\"older's inequality yields for the first term on the right-hand side
\begin{eqnarray*}
&&
T_n^{-1/2} a \sum_{i=0}^{n-1}E
\biggl[ \biggl\llvert X_{t_i} \int_{t_i}^{t_{i+1}}
(X_s - X_{t_i}) \dd s \biggr\rrvert\1_{\{|\di W + \di D|> v_n\}}
\biggr]
\\
&&\quad\leq T_n^{-1/2} a \sum_{i=0}^{n-1}E
\biggl[ \biggl( X_{t_i} \int_{t_i}^{t_{i+1}}
(X_s - X_{t_i}) \dd s \biggr)^2
\biggr]^{1/2} P\bigl(|\di W + \di D|> v_n\bigr)^{1/2}
\\
&&\quad= \mathrm{O}\bigl(T_n^{1/2} \dn^{1- \beta/2}\bigr)
\end{eqnarray*}
for the second addend we find that
\begin{eqnarray*}
&&
T_n^{-1/2} \dn a \sum_{i=0}^{n-1}E
\bigl[ X_{t_i}^2 \1_{\{|\di W + \di D|>
v_n\}}\bigr]
\\
&&\quad\leq T_n^{-1/2} \dn a \sum_{i=0}^{n-1}E
\bigl[ X_{t_i}^4\bigr]^{1/2} P\bigl(|\di W + \di D|>
v_n\bigr)^{1/2}
\\
&&\quad= \mathrm{O}\bigl(T_n^{1/2} \dn^{1/2 - \beta/2}\bigr).
\end{eqnarray*}
For $A_n^3$ we get by a similar estimate as for $A_n^1$ that
\[
E\bigl[\bigl|A_n^3\bigr|\bigr] = \mathrm{O} \bigl(T_n^{1/2}
\dn^{1/2 - \beta/2}\bigr).
\]
The last addend $A_n^4$ converges to zero in probability, since by
independence and H\"older's inequality
\begin{eqnarray*}
E\bigl[\bigl|A_n^4\bigr|\bigr] &\leq& T_n^{-1/2}
\sum_{i=0}^{n-1}E\bigl[X_{t_i}^2
\bigr]^{1/2} E\bigl[\di U^2\bigr]^{1/2} P\bigl(|\di W +
\di D|> v_n\bigr)^{1/2}
\\
&=& \mathrm{O}\bigl(T_n^{1/2} \dn^{1/2 - \beta/2}\bigr).
\end{eqnarray*}
\upqed\end{pf}

Now we show that the increments of the continuous part of $X$ are
negligible in the limit. This convergence is mainly based in the moment
bound that we have derived in Lemma \ref{pbound}.
%
\begin{lemma} \label{lemW}
\[
T_n^{-1/2} \sum_{i=0}^{n-1}X_{t_i}
\di X \1_{\{|\di W + \di D| > v_n\}
} \stackrel{p} {\longrightarrow}0 \qquad\mbox{as $n \to
\infty$}.
\]
\end{lemma}
\begin{pf}
We decompose $\di X = \di W + \di D + \di M + \di U$ to obtain
\begin{eqnarray*}
&&
T_n^{-1/2} \sum_{i=0}^{n-1}X_{t_i}
\di X \1_{\{|\di W + \di D| > v_n\}
} \\
&&\quad= T_n^{-1/2} \sum
_{i=0}^{n-1}X_{t_i} \di W \1_{\{|\di W + \di D| > v_n\}
}
+ T_n^{-1/2} \sum_{i=0}^{n-1}X_{t_i}
\di D \1_{\{|\di W + \di
D| > v_n\}} \\
&&\qquad{}+ T_n^{-1/2} \sum
_{i=0}^{n-1}X_{t_i} \di M \1_{\{|\di W + \di D| >
v_n\}}
+ T_n^{-1/2} \sum_{i=0}^{n-1}X_{t_i}
\di U \1_{\{|\di W + \di
D| > v_n\}}
\\
&&\quad= V_n^1 + V_n^2 +
V_n^3 + V_n^4.
\end{eqnarray*}
Lemma \ref{pbound} yields for $\delta= 1/2 - \beta$ and $l=2$ that
\[
P\bigl( |\di W + \di D| > v_n\bigr) = \mathrm{O}\bigl(\dn^{2- 2\beta}\bigr).
\]
For $V_n^1$ we obtain by H\"older's inequality and independence of
$X_{t_i}$ and $\di W$ that
\begin{eqnarray*}
E\bigl[\bigl|V_n^1\bigr|\bigr] &=& T_n^{-1/2}
E \Biggl[ \Biggl\llvert\sum_{i=0}^{n-1}X_{t_i}
\di W \1_{\{|\di
W + \di D| > v_n\}}\Biggr\rrvert\Biggr]
\\
&\leq& T_n^{-1/2} \dn^{1/2} \sum
_{i=0}^{n-1}E\bigl[X_{t_i}^2
\bigr]^{1/2} P\bigl(|\di W + \di D| > v_n\bigr)^{1/2}
\\
&=& \mathrm{O}\bigl(T_n^{1/2} \dn^{1/2 - \beta}\bigr).
\end{eqnarray*}
To prove convergence of $V_n^2$ we decompose $\di D$ as in (\ref
{Ddecomposition}) to obtain
\begin{eqnarray*}
E\bigl[\bigl|V_n^2\bigr|\bigr] &=& T_n^{-1/2}
E \Biggl[ \Biggl\llvert\sum_{i=0}^{n-1}X_{t_i}
\di D \1_{\{|\di
W + \di D| > v_n\}}\Biggr\rrvert\Biggr]
\\
&\leq& T_n^{-1/2}a E \Biggl[ \Biggl\llvert\sum
_{i=0}^{n-1}X_{t_i} \int
_{t_i}^{t_{i+1}} (X_s - X_{t_i})
\dd s \1_{\{|\di W + \di D| > v_n\}}\Biggr\rrvert\Biggr]
\\
&&{}+ T_n^{-1/2} a E \Biggl[ \Biggl\llvert\sum
_{i=0}^{n-1}X_{t_i}^2 \di
\1_{\{
|\di W + \di
D| > v_n\}}\Biggr\rrvert\Biggr].
\end{eqnarray*}
Applying H\"older's inequality to the first term on the right-hand side
results in
\begin{eqnarray*}
&&
T_n^{-1/2}a E \Biggl[ \Biggl\llvert\sum
_{i=0}^{n-1}X_{t_i} \int
_{t_i}^{t_{i+1}} (X_s - X_{t_i})
\dd s \1_{\{|\di W + \di D| > v_n\}}\Biggr\rrvert\Biggr]
\\
&&\quad\leq T_n^{-1/2}a \sum_{i=0}^{n-1}E
\biggl[ \biggl(X_{t_i} \int_{t_i}^{t_{i+1}}
(X_s - X_{t_i}) \dd s \biggr)^2
\biggr]^{1/2} P\bigl(|\di W + \di D| > v_n\bigr)^{1/2}
\\
&&\quad= \mathrm{O}\bigl(T_n^{1/2} \dn^{1-\beta}\bigr).
\end{eqnarray*}
The remaining term is of the order
\begin{eqnarray*}
&&
T_n^{-1/2} a E \Biggl[ \Biggl\llvert\sum
_{i=0}^{n-1}X_{t_i}^2 \di
\1_{\{|\di
W + \di
D| > v_n\}}\Biggr\rrvert\Biggr]
\\
&&\quad\leq T_n^{-1/2} a \dn\sum_{i=0}^{n-1}E
\bigl[ X_{t_i}^4 \bigr]^{1/2} P\bigl(|\di W + \di D| >
v_n\bigr)^{1/2}
= \mathrm{O}\bigl(T_n^{1/2} \dn^{1 - \beta}\bigr).
\end{eqnarray*}
Therefore, we conclude that $V_n^2 \to0$ as $n \to\infty$. Similar
estimates as for $V_n^1$ can be used for $V_n^3$ and $V_n^4$ to show
\[
E\bigl[\bigl|V_n^3\bigr|\bigr] = \mathrm{O}\bigl(T_n^{1/2}
\dn^{1- \beta}\bigr)
\quad\mbox{and}\quad
E\bigl[\bigl|V_n^4\bigr|\bigr] = \mathrm{O}\bigl(T_n^{1/2}
\dn^{1 -\beta}\bigr).
\]
This concludes the proof.
\end{pf}

The next lemma states that the increments of the jump component that
are close to the threshold are negligible in the limit. For the proof,
we use the small time moment bound for the jump component from
Proposition \ref{lemmoments}. This is the step where Assumption \ref
{aj} on the intensity of small jumps becomes crucial.
%
\begin{lemma} \label{lemM}
\[
T_n^{-1/2} \sum_{i=0}^{n-1}X_{t_i}
\di X \1_{\{ {v_n}/{2} < |\di
M| \leq2
v_n\}} \stackrel{p} {\longrightarrow}0 \qquad\mbox{as $n \to
\infty$}.
\]
\end{lemma}
\begin{pf}
Let us consider the following decomposition
\begin{eqnarray*}
&&
T_n^{-1/2} \sum_{i=0}^{n-1}X_{t_i}
\di X \1_{\{|\di M| > v_n/2, |\di
M| \leq2
v_n\}}
\\
&&\quad= T_n^{-1/2} \sum_{i=0}^{n-1}X_{t_i}
\di W \1_{\{|\di M| >
v_n/2, |\di
M| \leq2 v_n\}}
\\
&&\qquad{} + T_n^{-1/2} \sum_{i=0}^{n-1}X_{t_i}
\di D \1_{\{|\di M| >
v_n/2, |\di M|
\leq2 v_n\}}
\\
&&\qquad{} + T_n^{-1/2} \sum_{i=0}^{n-1}X_{t_i}
\di U \1_{\{|\di M| >
v_n/2, |\di M|
\leq2 v_n\}}
\\
&&\qquad{} + T_n^{-1/2} \sum_{i=0}^{n-1}X_{t_i}
\di M \1_{\{|\di M| >
v_n/2, |\di M|
\leq2 v_n\}}
\\
&&\quad= S_n^1 + S_n^2 +
S_n^3 + S_n^4.
\end{eqnarray*}
For the probability that $|\di M|$ lies in $(v_n/2, 2 v_n)$, we derive
from Proposition \ref{lemmoments} and Markov's inequality that
%
\begin{eqnarray}
\label{Mprob} P\bigl(|\di M| \leq2 v_n, | \di M| > v_n/2\bigr)
&=& P\bigl(|\di M| \1_{\{|\di M| \leq
2 v_n\}} > v_n/2\bigr)
\nonumber\\[-8pt]\\[-8pt]
&\leq& 4 v_n^{-2} E\bigl[(\di M)^2
\1_{\{|\di M| \leq2 v_n\}}\bigr]= \mathrm{O}\bigl(\dn^{1-\alpha\beta}\bigr).\nonumber
\end{eqnarray}
Hence, by independence of $X_{t_i}$, $\di W$, and $\di M$ we find that
$E[S_n^1] = 0$ and the second moment can be estimated as follows.
\begin{eqnarray*}
E\bigl[\bigl(S_n^1\bigr)^2\bigr] &\leq&
T_n^{-1} E \Biggl[\sum_{i=0}^{n-1}X_{t_i}^2
(\di W)^2 \1_{\{
|\di M| > v_n/2, |\di M| \leq2 v_n\}} \Biggr]
\\
&&{}+ T_n^{-1} E \biggl[ \sum_{ i \neq j}
X_{t_i} \di W \1_{\{|\di M| >
v_n/2, |\di M| \leq2 v_n\}} X_{t_j}
\Delta_j W \1_{\{|\Delta_j M| >
v_n/2, |\Delta_j M| \leq2 v_n\}} \biggr].
\end{eqnarray*}
Since $X_{t_i}, X_{t_j},\di W, \Delta_j M, \di M \perp\Delta_j W$
for $i < j$, the off-diagonal elements have zero expectation such that
the second addend vanishes. For the diagonal elements, we obtain
\begin{eqnarray*}
T_n^{-1} E \Biggl[\sum_{i=0}^{n-1}X_{t_i}^2
(\di W)^2 \1_{\{|\di M| >
v_n/2, |\di
M| \leq2 v_n\}} \Biggr] &\leq& T_n^{-1}
\dn\sum_{i=0}^{n-1}E\bigl[X_{t_i}^2
\bigr] \mathrm{O}\bigl(\dn^{1-\alpha\beta}\bigr)
\\
&=& \mathrm{O}\bigl(\dn^{1-\alpha\beta}\bigr).
\end{eqnarray*}
This yields the convergence $S_n^1 \stackrel{p}{\longrightarrow}0$ as
$n \to\infty$. To
prove that $S_n^2 \stackrel{p}{\longrightarrow}0$ as $n \to\infty$,
we plug in (\ref
{Ddecomposition}) and obtain
\begin{eqnarray*}
E\bigl[\bigl|S_n^2\bigr|\bigr] &\leq& E \Biggl[a
T_n^{-1/2} \Biggl\llvert\sum_{i=0}^{n-1}X_{t_i}
\int_{t_i}^{t_{i+1}} (X_s -
X_{t_i}) \dd s \1_{\{|\di M| > v_n/2, |\di M|
\leq2 v_n\}} \Biggr\rrvert\Biggr]
\\
&&{}+ E \Biggl[a T_n^{-1/2} \Biggl\llvert\sum
_{i=0}^{n-1}X_{t_i}^2 \di
\1_{\{
|\di M| >
v_n/2, |\di M| \leq2 v_n\}} \Biggr\rrvert\Biggr]
\end{eqnarray*}
and by independence
\begin{eqnarray*}
&&
E \Biggl[a T_n^{-1/2} \Biggl\llvert\sum
_{i=0}^{n-1}X_{t_i}^2 \di
\1_{\{|\di
M| >
v_n/2, |\di M| \leq2 v_n\}} \Biggr\rrvert\Biggr]
\\
&&\quad\leq a T_n^{-1/2} \sum_{i=0}^{n-1}E
\bigl[X_{t_i}^2\bigr] \di P\bigl(v_n/2 < |\di M|
\leq2 v_n\bigr)
\\
&&\quad= \mathrm{O}\bigl(T_n^{1/2} \dn^{1- \alpha\beta}\bigr).
\end{eqnarray*}
For the second term, H\"older's inequality yields
\[
E \Biggl[a T_n^{-1/2} \Biggl\llvert\sum
_{i=0}^{n-1}\int_{t_i}^{t_{i+1}}
(X_s - X_{t_i}) \dd s \1_{\{|\di M| > v_n/2, |\di M| \leq2 v_n\}} \Biggr
\rrvert
\Biggr] = \mathrm{O}\bigl(T_n^{1/2} \dn^{({1- \alpha\beta})/{2}}\bigr).
\]
From Assumption \ref{aj}, it follows that $S_n^4$ is centered for $n$
large enough. Furthermore, from Lemma \ref{lemj} we conclude
\begin{eqnarray*}
E\bigl[\bigl(S_n^4\bigr)^2\bigr] &=&
T_n^{-1} \sum_{i=0}^{n-1}E
\bigl[X_{t_i}^2\bigr] E\bigl[(\di M)^2 \1
_{\{|\di M| >
v_n/2, |\di M| \leq2 v_n\}}\bigr]
\\
&\leq& T_n^{-1} \sum_{i=0}^{n-1}E
\bigl[X_{t_i}^2\bigr] E\bigl[(\di M)^2
\1_{\{ |\di M|
\leq2
v_n\}}\bigr] \leq \mathrm{O}\bigl(\dn^{(2-\alpha) \beta}\bigr) \stackrel{n
\rightarrow\infty} {\longrightarrow}0.
\end{eqnarray*}
Finally, we show that $S_n^3 \stackrel{n\rightarrow\infty
}{\longrightarrow}0$. Independence together with
(\ref{Mprob}) leads to
\begin{eqnarray*}
E\bigl[\bigl|S_n^3\bigr|\bigr] &=& T_n^{-1/2}
\sum_{i=0}^{n-1}E\bigl[|X_{t_i} \di U
\1_{\{|\di
M| > v_n/2,
|\di M| \leq2 v_n\}}|\bigr]
\\
&=& T_n^{-1/2} \sum_{i=0}^{n-1}E\bigl[|X_{t_i}|\bigr]
E\bigl[|\di U|\bigr] P\bigl(|\di M| > v_n/2, |\di M| \leq2 v_n\bigr)
\\
&=& \mathrm{O}\bigl(T_n^{1/2} \dn^{1 - \alpha\beta}\bigr).
\end{eqnarray*}
\upqed\end{pf}
\begin{pf*}{Proof of Theorem \ref{discreteeff}}
Recall that for
\[
\hat{a}_n = - \frac{\sum_{i=1}^n X_{t_i} \Delta_i X^c}{\sum_{i=1}^n X_{t_i}^2
\Delta_i^n}
\]
we already know that $T_n^{1/2} (\hat{a}_n - a) \stackrel{\mathcal
{D}}{\longrightarrow}N(0, \frac
{\sigma^2}{E_a[X_\infty^2]})$ as $n \to\infty$. Therefore, it
remains to show
%
\begin{equation}
\label{aconv} T_n^{1/2} (\bar{a}_n -
\hat{a}_n) \stackrel{p} {\longrightarrow}0 \qquad\mbox{as $n \to\infty$.}
\end{equation}
Observe that
\[
T_n^{1/2} (\bar{a}_n - \hat{a}_n)
= \frac{T_n^{-1/2} (\sum_{i=1}^n
X_{t_i} \Delta_i X \1_{\{|\Delta_i X| \leq v_n\}} - \sum_{i=1}^n
X_{t_i} \Delta_i X^c)}{T_n^{-1} \sum_{i=1}^n X_{t_i}^2 \Delta_i}.
\]
By Lemma \ref{jumpfilter}
\[
T_n^{-1/2} \Biggl(\sum_{i=1}^n
X_{t_i} \Delta_i X \1_{\{|\Delta_i X|
\leq v_n\}} - \sum
_{i=1}^n X_{t_i} \Delta_i
X^c \Biggr) \stackrel{p} {\longrightarrow}0 \qquad\mbox{as $n\to\infty$}
\]
and
\[
T_n^{-1} \sum_{i=1}^n
X_{t_i}^2 \Delta_i^n \stackrel
{p} {\longrightarrow}E_a \bigl[X_\infty^2\bigr],
\]
such that (\ref{aconv}) follows.
\end{pf*}

\section{Simulation results} \label{chp7}

We investigate the finite sample performance of the estimator from
Sections \ref{chp5} and \ref{chp6} by means of Monte Carlo
simulations. First, we consider Ornstein--Uhlenbeck type processes with
finite jump intensity and give mean and standard deviation as well as
the number of jumps detected for different parameter values and varying
jump intensity. We also take a look at the normalized distribution of
the estimation error for finite samples. Then we investigate models
with infinite jump activity. In the last part, we compare the
performance of the maximum likelihood approach and least squares
estimation and find that
the jump filtering approach leads to a major improvement of the
estimate also for finite samples.

\subsection{Finite intensity models}

In this section, we perform Monte Carlo simulations for the drift
estimator (\ref{baranchp5}) of an Ornstein--Uhlenbeck type process
defined by
%
\begin{equation}
\label{OUsolchp7} X_t = \RMe^{-a t} X_0 + \int
_0^t \RMe^{-a (t-s)} \dd L_s,\qquad t \in
\R_+.
\end{equation}
We take a deterministic starting value $X_0 \in\R$ and $a > 0$. The
driving L\'evy process $L$ is assumed to be of the form
\[
L_t = W_t + \sum_{i=1}^{N_t}
Y_i,
\]
where $W$ is a Wiener process with $E[W_t^2] = \sigma_W^2 t$ and $N$
is a Poisson process with intensity $\lambda$ and the jump heights
$Y_i$ are i.i.d. with $N(0,2)$-distribution. An advantage of this
Ornstein--Uhlenbeck model is that exact simulation algorithms are
available both for $X$ and $L$. We use an exact discretization of the
explicit solution (\ref{OUsolchp7}) to the Langevin equation driven by
$L$ on a equidistant time grid $t_i = \dn i$ for $i = 1, \ldots, n$.
Algorithms for the exact simulation of $L$ can be found in \cite
{Cont}, among others.

Table \ref{tableOUPsim} contains means and standard deviations of each
$100$ realizations of the drift estimator $\bar a_n$ from (\ref
{baranchp5}). Since the Monte Carlo error is of order $N^{-1/2}$, where
$N$ is the number of Monte Carlo iterations, we have chosen a
reasonable compromise between precision of the Monte Carlo
approximation and computation time. The parameter values are $a=2$ and
$5$ and jump intensity~$\lambda$, time horizon $T$ and number of
observations $n$ vary as given in Table \ref{tableOUPsim}. We also
present the number of increments that were above the threshold $\dn
^{0.3}$. This number corresponds to the number of jumps that were
detected and we observe that it is relatively stable when $T$ and
$\lambda$ are kept fixed, which suggests that the jump filter works
quite reliable for finite intensity models and the threshold exponent
$\beta= 0.3$. For the compound Poisson process, the average number of
jumps in an interval of length $T$ is $E[N_T] = T \lambda$ and thus is
proportional
to the jump intensity. This relation is also visible for the simulated
data. The average number of filtered jumps is not equal to the expected
number of jumps, but lies between 60 and 70\% of the latter. This is
surprising, since we would expect the average number of detected jumps
to approach the expected number as $\dn$ tends to zero.

\begin{table}
\caption{Mean and standard deviation of $\bar a_n$ with $\beta= 0.3$
for an OU process with Gaussian component and compound Poisson
jumps}\label{tableOUPsim}
\begin{tabular*}{\tablewidth}{@{\extracolsep{\fill}}llllld{3.1}lld{3.1}@{}}
\hline
&&& \multicolumn{3}{l}{$a = 2$} & \multicolumn{3}{l@{}}{$a = 5$} \\[-4pt]
&&& \multicolumn{3}{l}{\hrulefill} & \multicolumn{3}{l@{}}{\hrulefill} \\
$\lambda$ & $T$ & $n$ & Mean & std dev & \multicolumn{1}{l}{$\varnothing$ jumps detect} &
Mean & std dev & \multicolumn{1}{l@{}}{$\varnothing$ jumps detect}\\
\hline
1 & 10 & 1000 & 2.0& 0.3 & 6.7& 5.0 &0.5 &7.4 \\
& & 2000 & 2.0& 0.3 & 7.0 & 5.0 &0.5 &7.2 \\
& & 4000 & 2.0 & 0.4 & 7.0& 5.0 &0.5 &6.8 \\
& 20 & 1000 & 2.0 & 0.2 & 13.1 & 4.7 & 0.3 & 12.5 \\
& & 2000 & 2.0 & 0.2 & 13.2 & 4.9 & 0.4 & 12.3 \\
& & 4000 & 2.0 & 0.2 & 13.0& 5.0 &0.3 & 13.1 \\
& 50 & 4000 & 2.0 & 0.1 & 31.3 &4.8 & 0.2 & 31.2 \\
& & 6000 & 2.0 & 0.2 & 32.2 & 4.6 & 0.3 & 30.1 \\
[6pt]
5& 10 & 1000 & 1.9 & 0.2 & 31.3 &4.6 &0.3 &30.0 \\
& & 2000 & 2.0 & 0.2 & 31.2 &4.8 &0.3 &30.9 \\
& & 4000 & 2.0 & 0.2 & 31.6 &4.9 &0.3 & 30.9 \\
& 20 & 2000 & 1.9 & 0.1 & 61.4 & 4.6 &0.2 & 60.2 \\
& & 4000 & 2.0 & 0.1 & 62.2 & 4.8 &0.2 &61.4 \\
& 50 & 4000 & 1.9 & 0.1 & 149 &4.6 &0.1 &145 \\
& & 6000 & 1.9 & 0.1 & 149 &4.7 &0.1 & 148 \\
\hline
\end{tabular*}
\end{table}

Another interesting finding is that as soon as the step size $\dn$ is
so small that the discretization error is negligible (cf. Section 4.2.4
in \cite{Mai2012} for an analysis of the discretization error), a
further increase in the number of observations does not improve the
accuracy of the estimator any further. This indicates that the
assumption of high-frequency observations is already reasonable when
the stochastic error dominates the discretization error at least for
finite intensity models.

\begin{figure}

\includegraphics{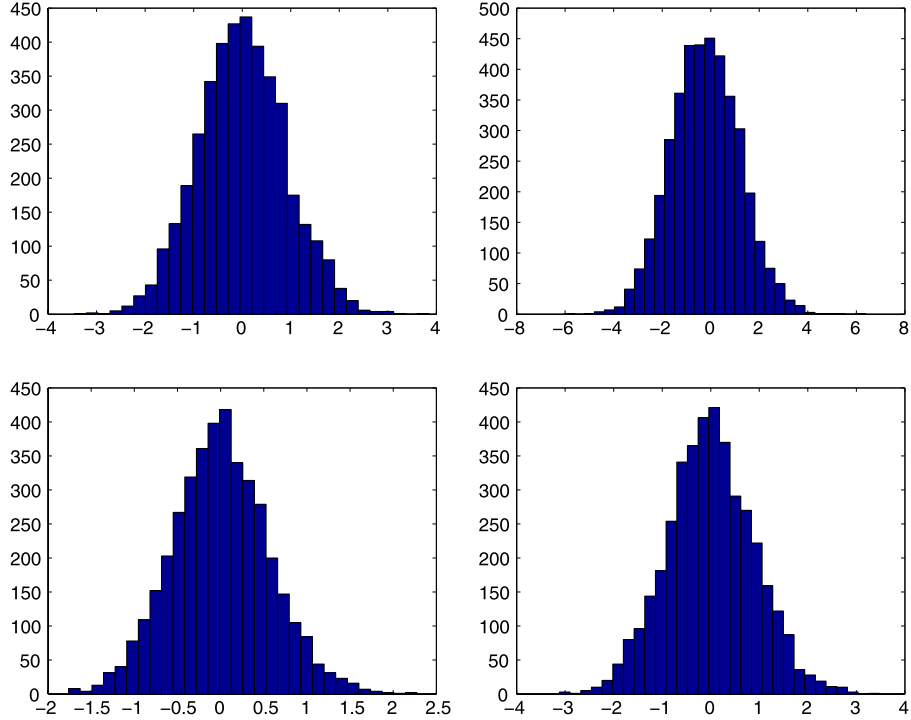}

\caption{Standardized error distribution of $\bar a_n$ for $a=2$
(left) and $a=5$ (right) compound Poisson jumps with intensity $\lambda
=1$ (top) and $\lambda= 2$ (bottom).} \label{OUPoissonhist1}
\end{figure}

The distribution of $T^{1/2} (\bar a_n - a)$ is depicted in Figure \ref
{OUPoissonhist1} for $T=70$, $\dn= 0.001$, $\sigma_W = 1$ and
$\lambda= 1$ and $2$. The histogram on the left corresponds to $a=2$
whereas on the right we have $a =5$. From Theorem \ref{discTCLT} and
the L\'evy--It\^o decomposition, it follows that the asymptotic variance
of $\bar a_n$ is given by
%
\begin{equation}
\label{avarsim1} \operatorname{AVAR} (\bar a_n) =\bigl(2a
\sigma_W^2\bigr) \bigl( \sigma_W^2
+ \lambda\sigma_j^2 \bigr)^{-1},
\end{equation}
where $\sigma_j^2$ denotes the variance of the jump heights. Hence, we
find that the asymptotic variance is proportional to $a$, which can
also be observed for finite samples in Figure \ref{OUPoissonhist1}.
By comparing the results of Figure \ref{OUPoissonhist1} (top) and
(bottom), we find that the variance also scales with the jump intensity
as indicated in (\ref{avarsim1}).

Eventually, we find that the estimator performs well even for very
short time horizons if the discretization is fine enough. This
observation corresponds to the results of Theorem 4.2.12 in \cite
{Mai2012} that states that the discretization bias is of the order
$\RMO(\dn)$.

\subsection{Infinite intensity models}

In Section \ref{chp6}, we have proved an asymptotic normality result
for the discretized maximum likelihood estimator with jump filter
(\ref{baranchp5}) for models that involve a jump component of infinite
activity. In this section, we simulate data from an Ornstein--Uhlenbeck
model of the form (\ref{OUsolchp7}) with $L = W + G$, where $W$ is a
Wiener process with $E [W_t^2] = \sigma_W^2 t$ and $G$ is a gamma
process. Again, we consider an equidistant grid $t_i = i \dn$ for $i =
0,\ldots, n$. The gamma process has jumps of infinite activity, paths
of finite variation and its Blumenthal--Getoor index is zero. The L\'evy
measure $\mu$ of $G$ has an explicit Lebesgue density given by
\[
g(x) = c x^{-1} \RMe^{-\lambda x} \1_{\{x >0\}} \qquad\mbox{for $x \in
\R$.}
\]
The parameter $c >0$ controls the jump intensity and $\lambda>0$ the
frequency of large jumps. It follows immediately from this density that
$G$ is a subordinator. Exact simulation algorithms are known for
increments of gamma processes and we use Johnk's algorithm (cf. \cite{Cont}).

\begin{figure}

\includegraphics{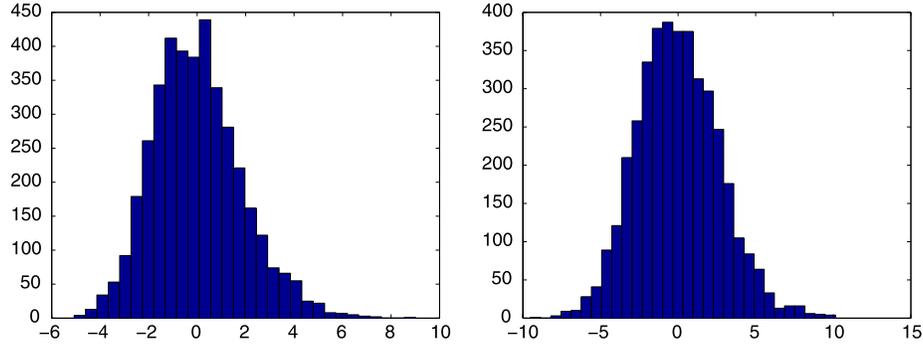}

\caption{Error distribution of $\bar a_n$ for an Ornstein--Uhlenbeck
process with $a=2$ (left) and $a=5$ (right), $\sigma_W = 1$ and gamma
process jumps.}\label{OUGammahist}
\end{figure}

Table \ref{ougammatable} gives mean and standard deviation for
different observation lengths and parameter values. The standard
deviation scales approximately with $T^{-1/2}$ as expected from Theorem
\ref{discreteeff}. In contrast to Table \ref{tableOUPsim}, we kept
here $\dn=0.0015$ fixed for all $n$. As in the finite intensity case
we use the threshold exponent $\beta= 0.3$ for the jump filter. We
find that the value of $a$ has hardly any impact on the average number
of increments that is filtered. When $a$ increases the number of
filtered increments also increases slightly, since a greater
variability of the drift might push increments with a relatively small
jump over the threshold.
Histograms of the standardized estimation error of $\bar {a}_n$ are given in Figure \ref{OUGammahist}
for $a = 2$ and $a = 5$ and jumps from a gamma process.

\begin{table}
\caption{Results of 200 Monte Carlo simulations of $\bar a_n$ with
$\Delta_n =0.0015$ and $\beta= 0.3$ for gamma process jumps}
\label{ougammatable}
\begin{tabular*}{\tablewidth}{@{\extracolsep{\fill}}ld{2.1}lld{2.1}lld{2.1}@{}}
\hline
&& \multicolumn{3}{l}{$a = 2$} & \multicolumn{3}{l@{}}{$a = 5$} \\[-4pt]
&& \multicolumn{3}{l}{\hrulefill} & \multicolumn{3}{l@{}}{\hrulefill} \\
$c$ & \multicolumn{1}{l}{$T$} & Mean & std dev & \multicolumn{1}{l}{$\varnothing$ jumps detect}
& Mean & std dev
& \multicolumn{1}{l@{}}{$\varnothing$ jumps detect}\\
\hline
0.5& 1 & 2.1 & 0.8 & 2.4 & 5.2 &1.2 & 2.2 \\
& 5 & 2.0 & 0.4 & 11.7 & 5.0 & 0.6 & 12.1 \\
& 7.5 & 2.0 & 0.3 & 17.7 & 4.9 & 0.5 & 17.8 \\
& 10 & 2.0 & 0.3 &23.7 & 5.0 & 0.4 &23.9 \\
& 20 & 2.0 & 0.2 & 47.2 & 5.0 & 0.3 & 47.6 \\
[6pt]
1& 1 & 2.1 & 0.8 & 1.6 &5.2 & 1.4& 1.8 \\
& 2.5 & 2.1 & 0.6 & 5.2 & 5.1 & 1.1 & 5.9 \\
& 5 & 2.1 & 0.5 & 8.4 & 5.0 & 0.8 & 8.6 \\
& 7.5 & 2.0 & 0.5 & 12.7 & 5.0 & 0.7 & 13.1 \\
& 10 & 2.0 & 0.3 & 17.2 & 5.0 & 0.6 & 17.1\\
\hline
\end{tabular*}
\end{table}

We conclude that the jump filtering approach performs well, also for
models with infinite jump activity provided that the maximal
observation distance is small.

\subsection{Maximum likelihood vs. least squares estimation} \label{sec73}

In this section, we compare maximum likelihood and least squares
estimation for the Ornstein--Uhlenbeck type process $X$ defined in
(\ref{OUsolchp7}). For continuously observed $X$, the least squares
estimator for the drift parameter $a$ is given by
\[
\hat{a}_T^{\mathrm{LS}} = -\frac{\int_0^T X_{s} \dd X_s}{\int_0^T X_{s}^2 \dd s}.
\]
For Gaussian Ornstein--Uhlenbeck processes, the least squares and the
likelihood estimator $\hat a_T^{\mathrm{ML}}$ (\ref{MLE}) coincide, since the
continuous martingale part under $P^a$ equals the process itself. But
when the driving process has jumps it follows from Theorem 4.2.10 in
\cite{Mai2012} that the asymptotic variances of both estimators differ by
\[
\operatorname{AVAR} \bigl(\hat a_T^{\mathrm{LS}} \bigr) -
\operatorname{AVAR} \bigl(\hat a_T^{\mathrm{ML}} \bigr) =
E_a \bigl[X_\infty^2\bigr]^{-1} \int
_\R x^2 \mu(\ddd x)> 0.
\]
Hence, the least squares estimator is inefficient when jumps are present.

%

\begin{figure}[b]

\includegraphics{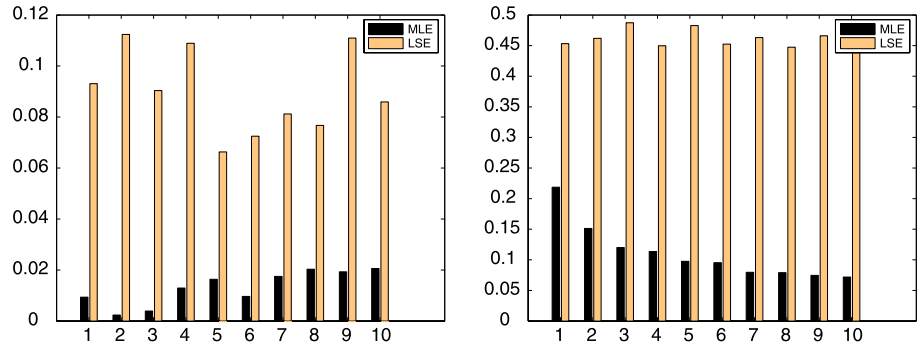}

\caption{Bias (left) and standard deviation (right) of MLE and LSE for
varying jump intensity.} \label{MLEvsLSEmean}
\end{figure}

Figure \ref{MLEvsLSEmean} compares the bias of the MLE and LSE for
compound Poisson jumps with varying jump intensity. For each intensity
the mean of 500 Monte Carlo simulations is given for an
Ornstein--Uhlenbeck process with $\sigma_W = 1$ and jumps with
$N(0,2)$-distribution. The true parameter is $a=2$ and we find that the
MLE has a significantly smaller bias than the LSE.

The standard deviation for both estimators is given in Figure \ref
{MLEvsLSEmean} on the right. The jump intensity $\lambda$ of the
compound Poisson part of $L$ varies between one and ten. In this model
setup, the difference in asymptotic variance between MLE and LSE is
given by
\[
\operatorname{AVAR} \bigl(\hat a_T^{\mathrm{LS}} \bigr) -
\operatorname{AVAR} \bigl(\hat a_T^{\mathrm{ML}} \bigr) =
\frac{2a \sigma_j^2 \lambda
}{\sigma_W^2 + \sigma_j^2 \lambda}.
\]
We observe that already for small jump intensities the MLE clearly
outperforms the LSE. With growing intensity this efficiency gain
becomes even more severe. For $\lambda= 10$, the standard deviation is
about five times larger for the least squares estimator.

This simulation example shows the significant gain in efficiency when
we use a discretized likelihood estimator with approximation of the
continuous part for drift estimation and underlines the importance of
jump filtering for jump diffusion models.

\section*{Acknowledgments}

The author is very grateful to Uwe K\"uchler and Markus Rei{\ss} for
stimulating comments and discussions and an anonymous referee for many
helpful questions and remarks that have led to considerable
improvements. This work was partially supported by Deutsche
Forschungsgemeinschaft through IRTG ``Stochastic Models of Complex
Processes''.



\printhistory

\end{document}